\documentclass{article}
\usepackage{amssymb}
\usepackage{graphicx}
\usepackage{amsmath}

\newtheorem{theorem}{Theorem}

\newtheorem{corollary}[theorem]{Corollary}

\newtheorem{definition}[theorem]{Definition}

\newtheorem{lemma}[theorem]{Lemma}

\newtheorem{proposition}[theorem]{Proposition}

\input{tcilatex}
\begin{document}

\title{Intrinsic characteristic classes\\
of a local Lie group \\
}
\author{Ender Abado\u{g}lu, Erc\"{u}ment Orta\c{c}gil}
\maketitle

\begin{abstract}
For a local Lie group $M$ we define cohomology classes $[w_{2k+1}]\in
H_{dR}^{2k+1}(M,\mathbb{R)}$. We show that $[w_{1}]$ is an obstruction to
globalizability and give an example where $[w_{1}]\neq 0.$ We also show that 
$[w_{3}]$ coincides with Godbillon-Vey class in a particular case. These
classes are secondary as they emerge when curvature vanishes.
\end{abstract}

\section{Introduction}

The problem studied in this paper emerged from a general framework which we
will outline in the next two paragraphs.

Let $M$ be a differentiable manifold with $\dim M=n\geq 2.$ Let $p,q\in M$
and $j_{k}(f)^{p,q}$ be the $k$-jet of a local diffeomorphism $f$ with $%
f(p)=q.$ We call $j_{k}(f)^{p,q}$ a $k$-arrow (from $p$ to $q$). A $0$-arrow
is an ordered pair $(p,q).$ Now let $G$ be a connected Lie group which acts
effectively on $G/H.$ Following a program outlined in [21], we defined in
[1] the geometric order $m$ of the (global) Klein geometry $G/H:$ $G/H$ has
order $m,$ if any $g\in G,$ as a transformation of $G/H,$ is determined by
any of its $m$-arrows and $m$ is the smallest such integer. Therefore $m=0$
eeeeeeeiff $H=\{e\}$. For any integer $m\geq 0,$ there exists a compact
Klein geometry $G/H$ of order $m.$ In fact, we can choose $G\subset SL(m,%
\mathbb{R})$ and $H=$ the Borel subgroup of upper triangular matrices ([1]).
If $\Delta \subset G$ is a discrete subgroup with $\Delta \cap H=\{e\},$ we
obtain a pseudogroup on the discrete quotient $\Delta \backslash G/H$ with
the property that its local diffeomorphisms are determined on their domains
by any of their $m$-arrows. In this way we obtain pseudogroups which are
finer than affine or projective pseudogroups. These pseudogroups exist even
on Riemann surfaces but they are always subordinate to an affine or
projective structure in this case.

Abstracting the arrows of the action of $G$ on $G/H,$ we introduced in [22]
the concept of a pre-homogeneous geometric structure of order $m$ on a
manifold $M$ and the curvature of such a structure. For $m=0,$ such a
structure is absolute parallelism studied in this paper. Riemannian
structures arise for $m=1$ by abstracting the $1$-arrows of the action of $%
SO(n)\propto \mathbb{R}^{n}$ on $(SO(n)\propto \mathbb{R}^{n})/SO(n)=\mathbb{%
R}^{n}.$ If $G/H$ has order $m,$ then $g$ is determined also by its $(m+1)$%
-arrows. This fact implies the existence of some canonical splittings and
associates some canonical \textquotedblleft torsionfree
connections\textquotedblright\ to such structures. As in this paper,
curvature vanishes if and only if $(m+1)$-arrows integrate to a pseudogroup
on $M$. If further the geometry is complete, then $M$ becomes a discrete
quotient of a global model $G/H$ (see Proposition 15). Therefore, curvature
is a measure of how much the geometry deviates from some local homogeneous
model. This way of looking at curvature parallels the one proposed in the
recent works [2], [3].

This paper is motivated by our attempt to do something globally interesting
with the curvature in [22] for $m=0.$ This simplest case incorporates all
the the main ideas and technical aspects of the above approach. The main
technical result of this paper is the construction of certain secondary
characteristic classes in the Lie algebra cohomology (therefore also in the
de Rham cohomology) of a local Lie group $M$. It turns out that the first
two classes are actually well known: $w_{1}$ is the character of the adjoint
representation and therefore $[w_{1}]$ is an obstruction to unimodularity.
We show that $[w_{1}]$ is also an obstruction to globalizability. The
condition $w_{3}=0$ is the well known Cartan's criterion for solvability and
the class $[w_{3}]$ is constructed first in [5] to show the nonvanishing of
the third cohomology group of a semisimple Lie algebra (see Proposition 29).
We find it surprising that $[w_{1}],$ $[w_{3}]$ have not been viewed so far,
to our knowledge, as the first two of a sequence of secondary characteristic
classes in the literature. On the other hand, it is shown in [20] that the
theory of local Lie groups is not a simple consequence of the global theory,
but has its own set of interesting and delicate geometric structures. In
this direction, we show here that it is possible to free the concept of a
local Lie group from being local in such a way that it contains the theory
of global Lie groups as a special case.

This paper is organized as follows. In Section 2 we recall the bracket $[$ , 
$]$ on sections of $J_{1}T=$ the algebroid of $\mathcal{U}_{1}$ as defined
in [17], [23], [24] and compute the exterior derivative of $1$-forms
explicitly in local coordinates.

In Section 3 we formulate the existence of a parallelism on $M$ as a
geometric structure defined by a splitting $\varepsilon $ of the groupoid
projection $\pi :\mathcal{U}_{1}\rightarrow \mathcal{U}_{0}.$

The splitting $\varepsilon :\mathcal{U}_{0}\rightarrow \mathcal{U}_{1}$
determines two splittings $\widetilde{\Gamma },$ $\widehat{\Gamma }%
:T\rightarrow J_{1}T$ on the level of algebroids. The pair $(\widetilde{%
\Gamma },\widehat{\Gamma })$ is a dual pair of affine connections as defined
in [2], [3]. In Section 4 we associate two curvatures $\mathcal{R}%
_{1}(\varepsilon ),$ $\mathcal{R}_{2}(\varepsilon )$ to $\varepsilon $ and
compute them explicitly in coordinates. We have $\mathcal{R}_{1}(\varepsilon
)=0$ but \textit{not necessarily} $\mathcal{R}_{2}(\varepsilon )=0.$

In Section 5 we construct a closed $1$-form $\omega _{1}$ in the algebroid
cohomology of $J_{1}T$ with trivial coefficients (see [18], [19], [6], [7],
[11] for cohomology of Lie algebroids). Our local computations uncover the
fundamental roles of the Spencer operator and the algebraic bracket in the
construction of this secondary characteristic class. It seems to us that the
construction of $\omega _{1}$ may be recast in the more general framework of
secondary characteristic classes for Lie algebroids introduced recently in
the works [6], [7], [11]. If this is the case, we hope that our approach
will have something to add to these works. However, our main concern here is
not $\omega _{1}$ but its restriction to $\widetilde{\Gamma }(T)\subset
J_{1}T$ which pulls back to a $1$-form $w$ in de Rham complex, not
necessarily closed. We prove $dw_{1}=Tr\mathcal{R}_{2}(\varepsilon )$
(Proposition 7). Our approach gives a totally new way of looking at torsion,
a point emphasized first in [28] in a somewhat different setting.

In Section 6 we show that $\mathcal{R}_{2}(\varepsilon )=0$ if and only if $%
M $ is a local Lie group (Proposition 8). Therefore the $1$-form $w_{1}$ is
closed on a local Lie group and the problem is to understand the obstruction
to its exactness.

In Section 7 we define another curvature $\mathcal{R}(\varepsilon )$ which
is the curvature mentioned in the second paragraph above. We show $\mathcal{R%
}(\varepsilon )=0$ $\Longleftrightarrow $ $\mathcal{R}_{2}(\varepsilon )=0.$
It turns out that the implication $\mathcal{R}_{2}(\varepsilon
)=0\Longrightarrow $ $\mathcal{R}(\varepsilon )=0$ is equivalent to the
traditional form of the third fundamental theorem of Lie. We construct an
explicit local primitive of the closed form $w_{1}$ (Proposition 21). This
fact implies the main result of this section (Proposition 22) that if $%
[w_{1}]\neq 0$ in $H_{dR}^{1}(M,\mathbb{R})$ then $M$ is not globalizable.
We give an example where $[w_{1}]\neq 0.$ Definitely, there is a relation
between the cohomology class $[w_{1}]$ and the modular class introduced in
[10] and studied further in [29]. As shown in [20], the concepts of
globalizability and global associativity are equivalent for local Lie
groups. Many explicit and nontrivial examples of local Lie groups are
constructed in [20] which fail to be globalizable. One of these examples is
even simply connected but this pathology can not occur if $M$ is complete in
view of Corollary 19. We plan to take up the study of these examples in some
future work.

In Section 8 we define the higher order analogs of the closed forms $\omega
_{1}$ and $w_{1}.$ If $G=SL(2,\mathbb{R)}$, $\Delta \subset G$ a cocompact
discrete subgroup and $\mathcal{F}$ is the standard codimension one
foliation on the local Lie group $G/\Delta ,$ we show (Proposition 28) that $%
[w_{3}]\in H^{3}(G/\Delta ,\mathbb{R)}$ coincides with the Godbillon-Vey
class $GV(\mathcal{F})$ and therefore $[w_{3}]\neq 0$ by the well known
Roussarie calculation (see [4]). We find this example particulary
interesting since our construction of $[w_{3}]$ is independent of any
foliation. It is natural to expect that these higher classes are also
obstructions to globalizability. We give an example, communicated to us by
the referee, which shows that this expectation is unjustified.

Therefore, generalizing the construction of the characteristic classes in
this paper to arbitrary geometric order $m,$ showing their nontriviality and
clarifying their geometric meaning remain as challenging problems.

\section{The algebroid $J_{1}T$}

Let $M$ be a (connected) differentiable manifold with $\dim M=n\geq 2.$ Let $%
p,q\in M$ and $j_{k}(f)^{p,q}$ be the $k$-jet of a local diffeomorphism $f$
with $f(p)=q.$ We call $j_{k}(f)^{p,q}$ a $k$-arrow (from $p$ to $q$). Note
that a $0$-arrow is an ordered pair $(p,q).$ Let $\mathcal{U}_{k}^{p,q}$
denote the set of all $k$-arrows on $M$ from $p$ to $q$. We define the set $%
\mathcal{U}_{k}\overset{def}{=}\cup _{p,q\in M}\mathcal{U}_{k}^{p,q}$. We
have the composition map $\mathcal{U}_{k}^{q,r}\times \mathcal{U}%
_{k}^{p,q}\rightarrow \mathcal{U}_{k}^{p,r}$ defined by $j_{k}(g)^{q,r}\circ
j_{k}(f)^{p,q}\overset{def}{=}j_{k}(g\circ f)^{p,r}.$ The differentiable
structure on $M$ induces a differentiable structure on $\mathcal{U}_{k}$ and 
$\mathcal{U}_{k}$ is a transitive Lie equation (in finite form), which is a
very special groupoid (see [17], [23], [24] for Lie equations in finite and
infinitesimal forms and [18], [19] and the references therein for Lie
groupoids and algebroids). Note that $\mathcal{U}_{0}$ is the pair groupoid $%
M\times M$. We have the projection map $\pi _{k,j}:\mathcal{U}%
_{k}\rightarrow \mathcal{U}_{j},$ $j\leq k+1,$ induced by the projection of
jets and $\pi _{k,j}$ is a morphism of groupoids, that is, it preserves
composition and inversion of arrows. In this paper we need only $\mathcal{U}%
_{0}$ and $\mathcal{U}_{1}.$

Recall that the algebroid of $\mathcal{U}_{0}$ is the tangent bundle $%
T\rightarrow M.$ We now recall the algebroid of $\mathcal{U}_{1}$ (see [23],
[24] for further details). Let $J_{1}T\rightarrow M$ be the vector bundle
whose fiber over $p\in M$ consists of $1$-jets of vector fields at $p.$ We
denote sections of $J_{1}T\rightarrow M$ by $X,Y$ and sections of $%
T\rightarrow M$ by $\xi ,\eta .$ We sometimes use the same notation for a
section and its value at a point. We also use the same notation $E$ for both
the total space of a vector bundle $E\rightarrow M$ and for the space of
global sections of $E\rightarrow M.$

Now a section $X$ of $J_{1}T\rightarrow M$ over $(U,x^{i})$ is of the form $%
(X^{i}(x),X_{j}^{i}(x))$ and the projection $\pi :J_{1}T\rightarrow T$ is
given by $(X^{i},X_{j}^{i})\rightarrow (X^{i}).$ We have the Spencer
operator $D:J_{2}T\rightarrow T^{\ast }\otimes J_{1}T$ locally given by $%
(X^{i},X_{j}^{i},X_{jk}^{i})\rightarrow (\partial
_{j}X^{i}-X_{j}^{i},\partial _{k}X_{j}^{i}-X_{kj}^{i}).$ We have the
algebraic bracket $\{$ $,$ $\}_{p}:(J_{2}T)_{p}\times
(J_{2}T)_{p}\rightarrow (J_{1}T)_{p}$ whose coordinate formula is obtained
by differentiating the usual formula for the bracket of two vector fields
twice, evaluating at $p$ and replacing derivatives with jet coordinates.
Finally, we have the Spencer bracket $[$ $,$ $]$ on sections of $%
J_{1}T\rightarrow T$ defined by

\begin{equation}
\lbrack X,Y]=\{\widetilde{X},\widetilde{Y}\}+i_{\pi X}D(\widetilde{Y}%
)-i_{\pi Y}D(\widetilde{X})
\end{equation}

In (1), $\widetilde{X},$ $\widetilde{Y}$ are arbitrary lifts of $X$ and $Y$
to sections of $J_{2}T$ and $i_{\pi X}$ denotes contraction with respect to
the vector field $\pi X.$ The bracket $[X,Y]$ does not depend on the lifts $%
\widetilde{X},\widetilde{Y}.$ If $X=(X^{i},X_{j}^{i})$ and $%
Y=(Y^{i},Y_{j}^{i}),$ using (1) we compute

\begin{eqnarray}
\lbrack X,Y]^{i} &=&X^{a}\frac{\partial Y^{i}}{\partial x^{a}}-Y^{a}\frac{%
\partial X^{i}}{\partial x^{a}} \\
\lbrack X,Y]_{j}^{i} &=&X_{j}^{a}Y_{a}^{i}-Y_{j}^{a}X_{a}^{i}+X^{a}\frac{%
\partial Y_{j}^{i}}{\partial x^{a}}-Y^{a}\frac{\partial X_{j}^{i}}{\partial
x^{a}}  \notag
\end{eqnarray}%
where we use summation convention in (2). With the bracket given by (2), $%
J_{1}T\rightarrow T$ becomes the Lie algebroid of the Lie groupoid $\mathcal{%
U}_{1}.$ Note that $\pi :J_{1}T\rightarrow T$ becomes a homomorphism of
algebroids in view of (2).

As shown in [7], [3], for any Lie algebroid $\mathcal{A},$ there is a unique
Lie algebroid structure on $J_{1}\mathcal{A}$ compatible with $\mathcal{A}.$
If $\mathcal{A}=T,$ then the bracket of $J_{1}\mathcal{A}$ \ as defined in
these works coincides with (2), even though the role of the Spencer operator
is not evident at first sight.

We use the same notation $[$ , $]$ for the brackets of $J_{1}T$ and $T.$ For
a vector field $\xi ,$ let $j_{1}\xi $ denote the first prolongation of $\xi
.$ In coordinates, if $\xi =(X^{i}),$ then $j_{1}\xi $ $=(X^{i},\frac{%
\partial Xi}{\partial x^{j}})\in J_{1}T.$ We have

\begin{equation}
j_{1}[\xi ,\eta ]=[j_{1}\xi ,j_{1}\eta ]
\end{equation}%
that is, $[$ , $]$ respects prolongation. (3) is easily checked using (2).

We now define a representation of $J_{1}T$ on $C^{\infty }(M)$ by $X(f)%
\overset{def}{=}\mathcal{L}_{\pi X}(f)$ where $\mathcal{L}_{\pi X}$ denotes
Lie derivative with respect to the vector field $\pi X,$ and consider the
cohomology of $J_{1}T$ with respect to this representation (see [18], [19],
[6], [7], [11] for cohomology of Lie algebroids). Up to Section 8, we will
be interested only in $1$-forms. Let $(J_{1}T)^{\ast }$ $\rightarrow M$ be
the dual bundle of $J_{1}T\rightarrow M.$

\begin{definition}
A $1$-form (of the algebroid $J_{1}T)$ is a (smooth) section of $%
(J_{1}T)^{\ast }$ $\rightarrow M.$
\end{definition}

So a $1$-form $\omega $ pairs linearly with sections of $J_{1}T$ to
functions on $M.$ The $1$-form $\omega $ is locally of the form $(\omega
_{i}(x),\omega _{j}^{i}(x))$ and the pairing is given by

\begin{equation}
\omega (X)\overset{def}{=}X^{a}(x)\omega _{a}(x)+X_{b}^{a}(x)\omega
_{a}^{b}(x)
\end{equation}

It is easy to derive the transformation laws of the components of $X$ and $%
\omega $ which we will not do here. These transformation laws show that $%
(X^{i},0)$ has no invariant meaning whereas $(0,X_{j}^{i})$ does and belongs
to $Ker(\pi )$. In contrast, $(\omega ^{i},0)$ does have invariant meaning
and represents an ordinary $1$-form as a cochain in the deRham complex
whereas $(0,\omega _{j}^{i})$ has no invariant meaning. Therefore $\omega $
reduces to an ordinary $1$-form eeeee if and only if $\omega _{j}^{i}=0.$

Now we want to compute the exterior derivative $\delta (\omega )$ of $\omega 
$ in coordinates. Differentiating (4) and contracting with the vector field $%
\pi X$, we get

\begin{eqnarray}
\mathcal{L}_{\pi X}\omega (Y) &=&X^{c}\frac{\partial \omega (Y)}{\partial
x^{c}}=X^{c}\frac{\partial Y^{a}}{\partial x^{c}}\omega _{a}+X^{c}Y^{a}\frac{%
\partial \omega _{a}}{\partial x^{c}} \\
&&+X^{c}\frac{\partial Y_{b}^{a}}{\partial x^{c}}\omega
_{a}^{b}+X^{c}Y_{b}^{a}\frac{\partial \omega _{a}^{b}}{\partial x^{c}} 
\notag
\end{eqnarray}

Interchanging $X$ and $Y$ in (5) and subtracting the resulting formula from
(5), we get $\mathcal{L}_{\pi X}\omega (Y)-\mathcal{L}_{\pi Y}\omega (X)=$

\begin{eqnarray}
&&(X^{c}Y^{a}-Y^{c}X^{a})\frac{\partial \omega _{a}}{\partial x^{c}}%
+(X^{c}Y_{b}^{a}-Y^{c}X_{b}^{a})\frac{\partial \omega _{a}^{b}}{\partial
x^{c}} \\
&&-(X_{b}^{c}Y_{c}^{a}-Y_{b}^{c}X_{c}^{a})\omega _{a}^{b}+(X^{c}\frac{%
\partial Y^{a}}{\partial x^{c}}-Y^{c}\frac{\partial X^{a}}{\partial x^{c}}%
)\omega _{a}  \notag \\
&&+(X_{b}^{c}Y_{c}^{a}-Y_{b}^{c}X_{c}^{a}+X^{c}\frac{\partial Y_{b}^{a}}{%
\partial x^{c}}-Y^{c}\frac{\partial X_{b}^{a}}{\partial x^{c}})\omega
_{a}^{b}  \notag
\end{eqnarray}%
where we added and subtracted $(X_{b}^{c}Y_{c}^{a}-Y_{b}^{c}X_{c}^{a})\omega
_{a}^{b}$ in the LHS of (6). In view of (2) and (4), the sum of the last two
terms in (6) is $\eta \lbrack X,Y]$ so that (6) becomes the well known
formula

\begin{equation}
\mathcal{L}_{\pi X}\omega (Y)-\mathcal{L}_{\pi Y}\omega (X)=\delta \omega
(X,Y)+\omega \lbrack X,Y]
\end{equation}%
where

\begin{eqnarray}
\delta \omega (X,Y) &=&(X^{c}Y^{a}-Y^{c}X^{a})\frac{\partial \omega _{a}}{%
\partial x^{c}}+(X^{c}Y_{b}^{a}-Y^{c}X_{b}^{a})\frac{\partial \omega _{a}^{b}%
}{\partial x^{c}} \\
&&-(X_{b}^{c}Y_{c}^{a}-Y_{b}^{c}X_{c}^{a})\omega _{a}^{b}  \notag
\end{eqnarray}

If $\omega _{j}^{i}=\delta _{j}^{i}$, note that the last term in (8)
vanishes since $Tr(ab)=Tr(ba)$ where $Tr$ denotes trace.

\section{Parallelism as splitting}

\begin{definition}
A splitting of $\pi _{1,0}:\mathcal{U}_{1}\rightarrow \mathcal{U}_{0}$ is a
morphism of groupoids $\varepsilon :\mathcal{U}_{0}\rightarrow \mathcal{U}%
_{1}$ such that $\pi _{1,0}\circ \varepsilon =id.$
\end{definition}

So a splitting $\varepsilon $ assigns to any ordered pair $(p,q)$ a $1$%
-arrow from $p$ to $q$ and this assignment preserves composition and
inversion of arrows. We will denote $\varepsilon (p,q)$ by $\varepsilon
^{p,q}.$ Further, $\varepsilon $ is differentiable and $\varepsilon (%
\mathcal{U}_{0})$ is an imbedded submanifold of $\mathcal{U}_{1}$. Note that 
$\varepsilon ^{p,p}=1$-arrow of the identity map at $p.$ Also, a $1$-arrow $%
f^{p,q}$ from $p$ to $q$ determines an isomorphism $(f^{p,q})_{\ast
}:T_{p}\rightarrow T_{q}$ and if two $1$-arrows determine the same
isomorphism, then they are identical. This observation supplies the proof of
the following simple

\begin{lemma}
The following are equivalent:

$i)$ There exists a splitting $\varepsilon :\mathcal{U}_{0}\rightarrow 
\mathcal{U}_{1}$

$ii)$ $M$ is parallelizable
\end{lemma}

Henceforth, we always assume that $M$ is parallelizable with a splitting $%
\varepsilon $ which we fix once and for all. The splitting $\varepsilon $ is
the analog of a torsionfree connection. Note that geometry (= parallelism)
and connection (= $\varepsilon )$ are identical objects in the present
framework.

Let $(U,x^{i})$, $(V,y^{i})$ be any two coordinate patches on $M.$ The
splitting $\varepsilon $ is of the form $\varepsilon _{j}^{i}(x,y)$ on $%
U\times V$ . We call the differentiable functions $\varepsilon _{j}^{i}(x,y)$
the components of $\varepsilon $ on $U\times V$. The isomorphism $%
(\varepsilon ^{x,y})_{\ast }$ is given by $X^{i}(x)\rightarrow \varepsilon
_{j}^{i}(x,y)X^{j}(x)=Y^{i}(y)$ in coordinates. We will repeatedly use the
local formulas $\varepsilon _{a}^{i}(y,z)\varepsilon
_{j}^{a}(x,y)=\varepsilon _{j}^{i}(x,z),$ $\varepsilon _{j}^{i}(x,x)=\delta
_{j}^{i}$ in the following sections.

\section{The curvatures $\mathcal{R}_{1}$($\protect\varepsilon $), $\mathcal{%
R}_{2}$($\protect\varepsilon $)}

Let $p\in (U,x^{i}).$ We choose $X(p)\in T_{p}$ arbitrarily and define a
vector field $\xi =(X^{i})$ on $(U,x^{i})$ by

\begin{equation}
X^{i}(x)\overset{def}{=}\varepsilon _{a}^{i}(p,x)X^{a}(p)
\end{equation}

\begin{definition}
A vector field $\vartheta $ is called $\varepsilon $-invariant if $%
(\varepsilon ^{p,q})_{\ast }\vartheta (p)=\vartheta (q),$ $p,q\in M.$
\end{definition}

We denote the Lie algebra of vector fields on $M$ by $\mathfrak{X}(M)$ and
the vector space of $\varepsilon $-invariant vector fields by $\mathfrak{X}%
_{\varepsilon }(M).$ Clearly $\xi \in \mathfrak{X}_{\varepsilon }(M)$ is
uniquely determined by $\xi (p)$ for any $p\in M$ and therefore $\dim 
\mathfrak{X}_{\varepsilon }(M)=\dim M.$ However, $\mathfrak{X}_{\varepsilon
}(M)$ need not be a Lie algebra, that is, we may not have $[\mathfrak{X}%
_{\varepsilon }(M),\mathfrak{X}_{\varepsilon }(M)]\subset \mathfrak{X}%
_{\varepsilon }(M).$ Some $\xi =(X^{i})\in \mathfrak{X}_{\varepsilon }(M)$
is given by (9) on $(U,x^{i}).$

Differentiating (9) and evaluating at $x=p,$ we obtain

\QTP{Body Math}
\begin{equation}
\frac{\partial X^{i}}{\partial x^{j}}(p)=\left[ \frac{\partial \varepsilon
_{a}^{i}(p,x)}{\partial x^{j}}\right] _{x=p}X^{a}(p)
\end{equation}%
We define

\QTP{Body Math}
\begin{equation}
\Gamma _{jk}^{i}(x)\overset{def}{=}\left[ \frac{\partial \varepsilon
_{k}^{i}(x,y)}{\partial y^{j}}\right] _{y=x}
\end{equation}

\QTP{Body Math}
Note that $\Gamma _{jk}^{i}(x)$ need not be symmetric in $j,k.$ This fact
will be of fundamental importance below. (10), (11) show that $\Gamma $
defines a map of vector bundles

\QTP{Body Math}
\begin{eqnarray}
\widetilde{\Gamma } &:&T\rightarrow J_{1}T \\
&:&X^{i}\longrightarrow (X^{i},\Gamma _{ja}^{i}X^{a})  \notag
\end{eqnarray}

\QTP{Body Math}
We have the exact sequence

\QTP{Body Math}
\begin{equation}
0\longrightarrow Ker(\pi )\longrightarrow J_{1}T\overset{\pi }{%
\longrightarrow }T\longrightarrow 0
\end{equation}%
where $Ker(\pi )\simeq T^{\ast }\otimes T.$ So $\pi \circ \widetilde{\Gamma }%
=id_{T}$ and $\widetilde{\Gamma }$ defines a splitting of the extension (13).

\QTP{Body Math}
Now let $\vartheta =(X^{i})$ be a vector field. We have the $PDE$ defined on 
$\mathfrak{X}(M)$ by

\QTP{Body Math}
\begin{eqnarray}
j_{1}(\vartheta ) &=&\widetilde{\Gamma }(\vartheta ) \\
\frac{\partial X^{i}}{\partial x^{j}} &=&\Gamma _{ja}^{i}X^{a}  \notag
\end{eqnarray}

\QTP{Body Math}
We define $\mathcal{R}_{1}(\varepsilon )$ by

\QTP{Body Math}
\begin{equation}
\mathcal{R}_{1}(\varepsilon )_{rj,k}^{i}\overset{def}{=}\left[ \frac{%
\partial \Gamma _{jk}^{i}}{\partial x^{r}}+\Gamma _{rk}^{a}\Gamma _{ja}^{i}%
\right] _{[r,j]}
\end{equation}%
where $[rj]$ denotes alternation. It is easy to show that $\mathcal{R}%
_{1}(\varepsilon )\in \wedge ^{2}(T^{\ast })\otimes T^{\ast }\otimes T.$ A
straightforward computation using (14) shows that $\mathcal{R}%
_{1}(\varepsilon )=0$ is the integrability condition of (14). Equivalently,
we can define $\nabla _{j}X^{i}\overset{def}{=}\frac{\partial X^{i}}{%
\partial x^{j}}-\Gamma _{ja}^{i}X^{a}\in T^{\ast }\otimes T$ which gives $%
\nabla _{r}\nabla _{j}X^{i}-\nabla _{j}\nabla _{r}X^{i}=\mathcal{R}%
_{1}(\varepsilon )_{rj,a}^{i}X^{a}.$ Note that the sign of $\Gamma
_{ja}^{i}X^{a}$ we use in $\nabla _{j}X^{i}$ is the opposite of the one in
tensor calculus: we could define $\Gamma _{jk}^{i}(x)$ by $\left[ \frac{%
\partial \varepsilon _{k}^{i}(x,y)}{\partial x^{j}}\right] _{y=x}$ .
Differentiating $\varepsilon _{a}^{i}(y,x)\varepsilon _{j}^{a}(x,y)=\delta
_{j}^{i}$ and evaluating on diagonal we get $\left[ \frac{\partial
\varepsilon _{k}^{i}(x,y)}{\partial y^{j}}\right] _{y=x}=-\left[ \frac{%
\partial \varepsilon _{k}^{i}(x,y)}{\partial x^{j}}\right] _{y=x}.$

\begin{lemma}
$\xi \in $ $\mathfrak{X}(M)$ belongs to $\mathfrak{X}_{\varepsilon }(M)$ if
and only if it is a solution of (14). In particular $\mathcal{R}%
_{1}(\varepsilon )=0.$
\end{lemma}

\QTP{Body Math}
Proof: Since $p$ is arbitrary in (10), any $\xi \in \mathfrak{X}%
_{\varepsilon }(M)$ is a solution of (14). Therefore $\mathcal{R}%
_{1}(\varepsilon )=0$ by the definition of $\mathcal{R}_{1}(\varepsilon ).$
Conversely, let $\vartheta \in \mathfrak{X}(M)$ be a solution of (14). We
choose $p\in (U,x^{i})$ and extend $\vartheta (p)$ to an $\varepsilon $%
-invariant $\xi $ on $M.$ Now $\vartheta ,\xi $ both solve (14) and have the
same initial condition at $p.$ By uniqueness, $\vartheta =\xi $ on some
neighborhood $p\in \overline{U}\subset U.$ Therefore, the set $A\overset{def}%
{=}\{p\in M\mid \vartheta (p)=\xi (p)\}\neq \emptyset $ is both open and
closed in $M$ and we conclude $A=M.$ \ $\square $

\QTP{Body Math}
To clarify the relation between (12), (14), (15) and the formalism of
connections on principal bundles, we now choose $p\in (U,x^{i})$ and define

\QTP{Body Math}
\begin{equation}
\overline{\Gamma }_{jk}^{i}(p,x)\overset{def}{=}\frac{\partial \varepsilon
_{b}^{i}(p,x)}{\partial x^{j}}\varepsilon _{k}^{b}(x,p)
\end{equation}

We claim that $\overline{\Gamma }_{jk}^{i}(p,x)$ is independent of $p.$
Differentiating $\varepsilon _{a}^{i}(q,x)\varepsilon
_{j}^{a}(p,q)=\varepsilon _{j}^{i}(p,x),$ we obtain

\QTP{Body Math}
\begin{equation}
\frac{\partial \varepsilon _{a}^{i}(q,x)}{\partial x^{k}}\varepsilon
_{j}^{a}(p,q)=\frac{\partial \varepsilon _{j}^{i}(p,x)}{\partial x^{k}}
\end{equation}

\QTP{Body Math}
Multiplying (17) with $\varepsilon _{b}^{i}(q,p)\varepsilon
_{j}^{b}(x,q)=\varepsilon _{j}^{i}(x,p)$ and summing over $j$ gives $%
\overline{\Gamma }_{jk}^{i}(q,x)=\overline{\Gamma }_{jk}^{i}(p,x)$ as
claimed. Now setting $p=x$ in (16), \ we conclude

\QTP{Body Math}
\begin{equation}
\overline{\Gamma }_{jk}^{i}(p,x)=\Gamma _{jk}^{i}(x)
\end{equation}

\QTP{Body Math}
\ \ We now fix $p\in M$ and a coordinate system around $p$ once and for all
and consider the principal bundle $\mathcal{U}_{1}^{(p)}\overset{def}{=}\cup
_{x\in M}\mathcal{U}^{p,x}$ $\rightarrow M$ with structure group $\mathcal{U}%
^{p,p}\simeq GL(n,\mathbb{R)}$. This principal bundle can be identified (not
canonically) with the principal frame bundle of $M.$ Now $\varepsilon $
trivializes this bundle as $x\rightarrow \varepsilon ^{p,x}$ and therefore
defines a flat connection with trivial monodromy. The components of this
connection are given by (18) and its curvature is $\mathcal{R}%
_{1}(\varepsilon ).$ However, the horizontal lift to the principal bundle $%
\mathcal{U}_{1}^{(p)}$ $\rightarrow M$ is \textit{not }given by (12) but by
(21) below. We define the Lie algebra bundle $\mathcal{L}\overset{def}{=}%
\cup _{q\in M}\mathcal{L}(\mathcal{U}^{q,q})$ where $\mathcal{L}(\mathcal{U}%
^{q,q})$ denotes the Lie algebra of $\mathcal{U}^{q,q}.$ Since $\mathcal{L}%
\simeq T^{\ast }\otimes T$ , we have $\mathcal{R}_{1}(\varepsilon )\in
\wedge ^{2}(T^{\ast })\otimes \mathcal{L}.$ This flat connection determines
a flat connection on the associated tangent bundle which is given by (12).
This gives another reason (more familiar than Lemma 5) why we should have $%
\mathcal{R}_{1}(\varepsilon )=0.$

\QTP{Body Math}
\ \ Now 
\begin{equation}
\Gamma _{jk}^{i}-\Gamma _{kj}^{i}\overset{def}{=}T_{jk}^{i}=\text{torsion}
\end{equation}

\QTP{Body Math}
As for $\mathcal{R}_{2}(\varepsilon ),$ following [21] (see also [3] for a
different but equivalent definition) we define

\QTP{Body Math}
\begin{equation}
\mathcal{L}_{X}\xi \overset{def}{=}[\pi X,\xi ]+i_{\xi }D(X)
\end{equation}%
where $D:J_{1}T\rightarrow T^{\ast }\otimes T$ is the Spencer operator
locally given by $(X^{i},X_{j}^{i})\rightarrow (\frac{\partial X^{i}}{%
\partial x^{j}}-X_{j}^{i})$. Therefore $\mathcal{L}_{X}:T\rightarrow T$ is a
first order differential operator. To emphasize the analogy with $\mathcal{R}%
_{1}(\varepsilon ),$ we denote $\mathcal{L}_{\Gamma (\frac{\partial }{%
\partial x^{i}})}$ by $\widehat{\nabla }_{i}.$ Using (2) and (20) we find $%
\widehat{\nabla }_{j}X^{i}=\frac{\partial X^{i}}{\partial x^{j}}-\Gamma
_{aj}^{i}X^{a}$ which gives another splitting $\widehat{\Gamma }$ of (13)
defined by

\QTP{Body Math}
\begin{eqnarray}
\widehat{\Gamma } &:&T\rightarrow J_{1}T \\
&:&X^{i}\longrightarrow (X^{i},\Gamma _{aj}^{i}X^{a})  \notag
\end{eqnarray}

\QTP{Body Math}
\ \ If $T_{jk}^{i}=0,$ then clearly $\widetilde{\Gamma }=\widehat{\Gamma }.$
We define

\QTP{Body Math}
\begin{equation}
\mathcal{R}_{2}(\varepsilon )_{rj,k}^{i}\overset{def}{=}\left[ \frac{%
\partial \Gamma _{kj}^{i}}{\partial x^{r}}+\Gamma _{kr}^{a}\Gamma _{aj}^{i}%
\right] _{[r,j]}
\end{equation}%
and check that $\mathcal{R}_{2}(\varepsilon )=0$ is the integrability
condition of $\frac{\partial X^{i}}{\partial x^{j}}=\Gamma _{aj}^{i}X^{a}.$
It is also straightforward to check $\widetilde{\nabla }_{r}\widetilde{%
\nabla }_{j}X^{k}-\widetilde{\nabla }_{j}\widetilde{\nabla }_{r}X^{k}=%
\mathcal{R}_{2}(\varepsilon )_{rj,a}^{k}X^{a}$. Like $\mathcal{R}%
_{1}(\varepsilon ),$ $\mathcal{R}_{2}(\varepsilon )\in \wedge ^{2}(T)\otimes
T^{\ast }\otimes T.$

\QTP{Body Math}
Clearly $\widetilde{\Gamma }[\xi ,\eta ]-[\widetilde{\Gamma }\xi ,\widetilde{%
\Gamma }\eta ]\in Ker(\pi )$ for $\xi ,\eta \in \mathfrak{X}(M).$ The same
statement holds for $\widehat{\Gamma }.$ Now a straightforward computation
using (2), (12) and (21) gives the formulas

\QTP{Body Math}
\begin{eqnarray}
(\widetilde{\Gamma }[\xi ,\eta ]-[\widetilde{\Gamma }\xi ,\widetilde{\Gamma }%
\eta ])_{j}^{i} &=&\mathcal{R}_{2}(\varepsilon )_{ab,j}^{i}X^{a}Y^{b}\qquad
\\
(\widehat{\Gamma }[\xi ,\eta ]-[\widehat{\Gamma }\xi ,\widehat{\Gamma }\eta
])_{j}^{i} &=&\mathcal{R}_{1}(\varepsilon )_{ab,j}^{i}X^{a}Y^{b}=0\text{ \ \
\ }
\end{eqnarray}%
where $\xi =(X^{i}),\eta =(Y^{i}).$

\QTP{Body Math}
\ \ Now (24) shows that $\widehat{\Gamma }(T)\subset J_{1}T$ is a
subalgebroid which is the algebroid of $\varepsilon (\mathcal{U}_{0})\subset 
\mathcal{U}_{1}.$ In fact, the splitting $\xi \rightarrow \widehat{\Gamma }%
\xi $ is the horizontal lift to the principal bundle $\mathcal{U}_{1}^{(p)}$ 
$\rightarrow M$ which is trivialized by $\varepsilon .$ On the other hand,
the splitting (12) is conceptually different as we will see in Sections 5,
6. The pair $(\widetilde{\Gamma },\widehat{\Gamma })$ is a dual pair of
affine connections according to [2], [3]. Clearly, $T_{jk}^{i}=0$ implies $%
\mathcal{R}_{2}(\varepsilon )=\mathcal{R}_{1}(\varepsilon )=0,$ but $%
\mathcal{R}_{2}(\varepsilon )$ need not vanish in general as we will see
below.

\QTP{Body Math}
Finally, using (2) and (20) it is easy to check

\QTP{Body Math}
\begin{equation}
\mathcal{L}_{X}\mathcal{L}_{Y}\xi -\mathcal{L}_{Y}\mathcal{L}_{X}\xi =%
\mathcal{L}_{[X,Y]}\xi
\end{equation}

\QTP{Body Math}
\ \ Therefore (25) defines a representation of $J_{1}T$ \ on $T$ . More
generally, $J_{k+1}T$ has a representation on $J_{k}T$ (called association
in [23], see Definition 8.1 on pg. 362) defined by $\mathcal{L}_{X}\xi 
\overset{def}{=}[\pi _{1}X,\xi ]-i_{\pi _{2}\xi }D(X)$ where $X,\xi $ are
sections of $J_{k+1}T\rightarrow M,$ $J_{k}T\rightarrow M,$ and $\pi
_{1}:J_{k+1}T\rightarrow J_{k}T,$ $\pi _{2}:J_{k+1}T\rightarrow T$ are
projections ([23], Lemma 8.32). This association is used in [23] to study
deformation cohomology and rigidity as the culmination of this book (pg.
354-393). These concepts are studied also in [8] for a general Lie algebroid 
$\mathcal{A},$ but it seems to us that the theory introduced in [8] reduces
to [23] if $\mathcal{A}\subset J_{k}T.$

\section{A closed 1-form}

\QTP{Body Math}
We define

\QTP{Body Math}
\begin{equation}
Tr\mathcal{R}_{2}(\varepsilon )\overset{def}{=}\mathcal{R}_{2}(\varepsilon
)_{rj,a}^{a}=\left[ \frac{\partial \Gamma _{aj}^{a}}{\partial x^{r}}\right]
_{[r,j]}
\end{equation}%
The second equality in (26) follows from (22) since $\left[ \Gamma
_{br}^{a}\Gamma _{aj}^{b}\right] _{[r,j]}=0.$

\QTP{Body Math}
We will now define a special $1$-form $\omega $ of $J_{1}T.$ Recall the
algebraic bracket

\QTP{Body Math}
\begin{equation}
\{\text{ },\text{ }\}_{p}:(J_{1}T)_{p}\times (J_{1}T)_{p}\rightarrow T_{p}
\end{equation}

\QTP{Body Math}
A section $X$ of $J_{1}T$ defines the linear map%
\begin{eqnarray}
X(p) &:&T_{p}\rightarrow T_{p} \\
&:&\xi (p)\rightarrow \{\text{ }X(p),\text{ }\widetilde{\Gamma }(\xi
(p))\}_{p}  \notag
\end{eqnarray}

\QTP{Body Math}
We define $\omega $ by

\QTP{Body Math}
\begin{equation}
\omega (X)\overset{def}{=}Tr(X)
\end{equation}

\QTP{Body Math}
Clearly $\omega $ is a $1$-form of $J_{1}T.$

\begin{proposition}
$\omega $ is closed but not locally exact.
\end{proposition}

\QTP{Body Math}
Proof: If $X=$ $(X^{i},X_{j}^{i}),$ $Y=(Y^{i},Y_{j}^{i}),$ then (27) is
given by

\QTP{Body Math}
\begin{equation}
\{X,Y\}^{i}=X^{a}Y_{a}^{i}-Y^{a}X_{a}^{i}
\end{equation}

\QTP{Body Math}
Thus we get $\{\xi ,$ $\widetilde{\Gamma }\eta \}^{i}=X^{a}\Gamma
_{ab}^{i}Y^{b}-Y^{a}X_{a}^{i}$ where $\xi =(X^{i}),$ $\eta =(Y^{i}).$ So the
the linear map (28) is $(X^{i},X_{j}^{i}):$ $Y^{i}\rightarrow X^{a}\Gamma
_{ab}^{i}Y^{b}-Y^{a}X_{a}^{i}$ and we deduce

\QTP{Body Math}
\begin{equation}
Tr(X)=X^{a}\Gamma _{ab}^{b}-X_{a}^{a}
\end{equation}

\QTP{Body Math}
$\bigskip $(4), (29) and (31) show%
\begin{equation}
\omega =(\omega ^{i},\omega _{j}^{i})=(\Gamma _{ia}^{a},-\delta _{j}^{i})
\end{equation}

\QTP{Body Math}
$\bigskip $Substituting (32) into (11), we get%
\begin{eqnarray}
\delta \omega (X,Y) &=&X^{a}Y^{c}(\frac{\partial \Gamma _{ab}^{b}}{\partial
x^{c}}-\frac{\partial \Gamma _{cb}^{b}}{\partial x^{a}})  \notag \\
&=&Tr\mathcal{R}_{1}(\varepsilon )(\pi X,\pi Y)  \notag \\
&=&0\text{ }
\end{eqnarray}

\QTP{Body Math}
If $\delta f=\omega $ locally$,$ (32) implies $(\frac{\partial f}{\partial
x^{i}},0)=(\Gamma _{ia}^{a},-\delta _{j}^{i})$, an equality which does not
hold for any $f.$ \ $\square $

\QTP{Body Math}
Proposition 6 shows that the complex computing the cohomology of an
algebroid need not be locally exact at the level of 1-forms, which, we
believe, is a serious defect from the point of view of characteristic
classes since a closed form in such a complex has no interpretation as an
obstruction. This fact makes $\omega $ uninteresting from our standpoint and
we will use it as a tool for a more relevant construction. For this purpose,
note that the algebraic bracket (27) defines the alternating bilinear map

\QTP{Body Math}
\begin{eqnarray}
\lbrack \text{ , }]_{p} &:&T_{p}\times T_{p}\longrightarrow T_{p} \\
&:&(\xi (p),\eta (p))\rightarrow \{\widetilde{\Gamma }(\xi (p)),\widetilde{%
\Gamma }(\eta (p))\}_{p}  \notag
\end{eqnarray}

\QTP{Body Math}
We define an ordinary $1$-form $w$ on $M$ by

\QTP{Body Math}
\begin{equation}
w(\xi )\overset{def}{=}Tr[\xi ,\bullet ]
\end{equation}

\QTP{Body Math}
Clearly, $w(\xi )=\omega (\widetilde{\Gamma }\xi ).$

\begin{proposition}
$dw=Tr\mathcal{R}_{2}(\varepsilon )$
\end{proposition}

\QTP{Body Math}
Proof: (35), (34) and (30) give $w(\xi )=(\Gamma _{ab}^{b}-\Gamma
_{ba}^{b})X^{a}$ or

\QTP{Body Math}
\begin{equation}
w_{i}=\Gamma _{ia}^{a}-\Gamma _{ai}^{a}
\end{equation}%
Using (36), we compute

\QTP{Body Math}
\begin{eqnarray}
(dw)_{ji} &=&\left[ \frac{\partial \Gamma _{ia}^{a}}{\partial x^{j}}-\frac{%
\partial \Gamma _{ai}^{a}}{\partial x^{j}}\right] _{[j,i]} \\
&=&\left[ \frac{\partial \Gamma _{ai}^{a}}{\partial x^{j}}\right] _{[j,i]}%
\text{ \ \ \ since }\left[ \frac{\partial \Gamma _{ia}^{a}}{\partial x^{j}}%
\right] _{[j,i]}=Tr(\mathcal{R}_{1}(\varepsilon ))_{ji}=0  \notag \\
&=&Tr\mathcal{R}_{2}(\varepsilon )_{ji}\text{ \ \ \ \ by (26) }\qquad \square
\notag
\end{eqnarray}%
The derivation of (36) shows that $T_{jk}^{i}(p)=\Gamma _{jk}^{i}(p)-\Gamma
_{kj}^{i}(p)$ is a $Hom(T_{p},T_{p})$-valued $1$-form in the present
framework: $j$ is the $1$-form index and $i,k$ are the matrix indices. This
interpretation of torsion will be of fundamental importance in the
construction of secondary characteristic classes in Section 8.

\QTP{Body Math}
Therefore, $w$ is closed in the deRham complex if and only if $Tr\mathcal{R}%
_{2}(\varepsilon )=0.$ If we assume $T_{jk}^{i}=0,$ then $dw=Tr\mathcal{R}%
_{2}(\varepsilon )=Tr\mathcal{R}_{1}(\varepsilon )=0.$ In fact, (19) and
(36) show that $T_{jk}^{i}=0$ implies $w=0$ $!$ The condition $T_{jk}^{i}=0$
is very strong: If $M$ is compact, it forces $M$ to be of the form $G/\Delta 
$ where $G$ is the affine linear group $AL(n,\mathbb{R)}$ and $\Delta
\subset G$ is a discrete subgroup (see [16], \ Theorem 4.2). However, our
assumption in the next section will allow $G$ to be \textit{any }connected
Lie group.

\section{From Lie algebroid to Lie algebra}

\begin{proposition}
The following are equivalent:

$i)$ $\mathcal{R}_{2}(\varepsilon )=0$

$ii)$ $\widetilde{\Gamma }(T)\subset J_{1}T$ is a subalgebroid

$iii)$ $[\mathfrak{X}_{\varepsilon }(M),\mathfrak{X}_{\varepsilon
}(M)]\subset \mathfrak{X}_{\varepsilon }(M)$
\end{proposition}

Proof: The equivalence $i)\Leftrightarrow ii)$ follows from (23).

$iii)\Rightarrow i):$ Let $\xi ,\eta \in \mathfrak{X}_{\varepsilon }(M).$ We
have $j_{1}[\vartheta ,\gamma ]=[j_{1}\vartheta ,j_{1}\gamma ]$ for any $%
\vartheta ,\gamma \in \mathfrak{X}(M)$ by (3). By Lemma 5, $j_{1}\vartheta =%
\widetilde{\Gamma }(\vartheta )$ for any $\vartheta \in \mathfrak{X}%
_{\varepsilon }(M)$ and therefore $j_{1}[\xi ,\eta ]=\widetilde{\Gamma }%
([\xi ,\eta ]$ since $[\xi ,\eta ]\in \mathfrak{X}_{\varepsilon }(M).$ So we
deduce $\widetilde{\Gamma }([\xi ,\eta ])=[\widetilde{\Gamma }\xi ,%
\widetilde{\Gamma }\eta ]$ \ $\xi ,\eta \in \mathfrak{X}_{\varepsilon }(M).$
If $(\xi _{i})$ is a basis for $\mathfrak{X}_{\varepsilon }(M),$ any $%
\vartheta \in \mathfrak{X}(M)$ can be written as $\vartheta =f^{i}\xi _{i}$
for some functions $f^{i}$ on $M$ and the conclusion follows since $\mathcal{%
R}_{2}(\varepsilon )(\vartheta ,\gamma )$ is linear in its arguments $%
\vartheta ,\gamma $.

$i)\Rightarrow iii):$ Let $\xi ,\eta \in \mathfrak{X}_{\varepsilon }(M)$.
(3), (23) and Lemma 5 give $j_{1}[\xi ,\eta ]=\widetilde{\Gamma }([\xi ,\eta
])$. Therefore $[\xi ,\eta ]$ is a solution of (14) and the conclusion
follows from Lemma 5 \ $\square $

We will denote the conditions of Proposition 8 by $\mathbf{A.}$ Propositions
7, 8 now give

\begin{corollary}
$\mathbf{A}$ implies $dw=0$
\end{corollary}

If $\xi _{(i)}$ is a basis of $\mathfrak{X}_{\varepsilon }(M)$, then $[\xi
_{(i)},\xi _{(j)}]^{k}=c_{ij}^{a}\xi _{(a)}^{k}$ for some \textit{functions }%
$c_{ij}^{a}$ on $M.$ These functions are constant if and only if $iii)$ of
Proposition 8 holds. We now make the following

\begin{definition}
If a differentiable manifold $M$ admits a splitting $\varepsilon $ such that 
$\mathbf{A}$ holds, then $M$ together with $\varepsilon $ is called a local
Lie group.
\end{definition}

\QTP{Body Math}
Our definition of local Lie group coincides with the one given in [20] in
view of Theorem 18 in [20] and Proposition 8. For a local Lie group $M$, we
will construct in Section 7 a Lie group $\widetilde{G}$ (see the proof of
Proposition 15) whose Lie algebra can be identified with $\mathfrak{X}%
_{\varepsilon }(M).$

\QTP{Body Math}
Henceforth in this section we assume that $M$ is a local Lie group.
Therefore $[w]\in H_{dR}^{1}(M,\mathbb{R)}$. Recalling the formula $w(\xi
)=\omega (\widetilde{\Gamma }\xi )$, we may identify $w$ with the
restriction of $\omega $ to $\widetilde{\Gamma }(T).$ Unlike $\omega ,$ $w$
is clearly locally exact and can be also globally exact in the deRham
complex. Our main concern in Section 7 will be the following question

$\mathbf{Q:}$ What is the obstruction to the exactness of $w$ $?$

We now define the evaluation map

\begin{eqnarray}
e_{p} &:&\mathfrak{X}_{\varepsilon }(M)\longrightarrow T_{p} \\
&:&\xi \longrightarrow \xi (p)  \notag
\end{eqnarray}%
where $p\in M$ is arbitrary. Recalling the alternating bilinear map $[$ , $%
]_{p}$ defined by (34), we have

\begin{proposition}
$e$ preserves brackets.
\end{proposition}

Proof: This is a straightforward verification using (9) and (34). \ $\square 
$

Proposition 11 shows that the Lie algebra structure of $\mathfrak{X}%
_{\varepsilon }(M)$ is determined at any point $p\in M$ in terms of $[$ , $%
]_{p}.$ The Jacobi identity imposes now restrictions on the components $%
\Gamma _{jk}^{i}(p)$ which we do not write down explicitly as we do not use
them. It is not difficult to check that we now have $(\varepsilon
^{p,q})_{\ast }[\xi ,\eta ]_{p}=[(\varepsilon ^{p,q})_{\ast }\xi
,(\varepsilon ^{p,q})_{\ast }\eta ]_{q}$ and $w$ is $\varepsilon $%
-invariant: $w_{p}(\xi (p))=w_{q}((\varepsilon ^{p,q})_{\ast }\xi (p)).$

\section{$[w_{1}]$ as obstruction to globalizability}

When is a manifold a Lie group? This question is quite old and studied by
several authors, for instance, see [12], [13], [26] for proofs based on the
formalism of Maurer-Cartan forms. A different proof is outlined in [2] which
is more in the spirit of this paper: assuming $\mathbf{A}$, the dual pair of
flat affine connections $(\widetilde{\Gamma },\widehat{\Gamma })$ imposes a 
\textit{local }Lie group structure on $M$ by Theorem A in [2]. We will give
a slightly different proof of this fact in this section which unifies the
dual pair $(\widetilde{\Gamma },\widehat{\Gamma })$ as the single object $%
\varepsilon $ and the curvatures $\mathcal{R}_{1}(\varepsilon ),$ $\mathcal{R%
}_{2}(\varepsilon )$ as $\mathcal{R}(\varepsilon )$ which is \textit{the
curvature }of $\varepsilon $ from our standpoint. Our proof will also
interpret the obstruction $[w]$ to globalizability of a local Lie group as a
secondary characteristic class in view of Corollary 9 and Proposition 7.

By the definition of a $1$-arrow, for any $p\in (U,x^{i})$ and $q\in
(V,y^{i})$, there exists a local diffeomorphism $f$ with $f(p)=q$ satisfying

\begin{equation}
\frac{\partial f^{i}}{\partial x^{j}}(p)=\varepsilon _{j}^{i}(p,f(p))
\end{equation}

Thus the splitting $\varepsilon $ determines a global $PDE$ on $\mathcal{U}%
_{0}=M\times M$ which is locally given by (39). Note that $f$ in (39)
depends on the point $(p,q)$ and we may not be able to find some $f$ which
works for all $p\in U$ as the $PDE$ (39) may not admit any local solutions.
The integrability condition of (39) is given by

\begin{equation}
\left[ \frac{\partial \varepsilon _{j}^{i}(x,y)}{\partial x^{k}}+\frac{%
\partial \varepsilon _{j}^{i}(x,y)}{\partial y^{a}}\varepsilon _{k}^{a}(x,y)%
\right] _{[kj]}\overset{def}{=}\mathcal{R}(\varepsilon )_{kj}^{i}(x,y)=0
\end{equation}%
If (39) admits a solution $f$ with $f(p)=q$ for any $(p,q)\in U\times V$,
then clearly $\mathcal{R}(\varepsilon )_{kj}^{i}(x,y)=0$ on $U\times V.$
Conversely, by the well known existence and uniqueness theorem for first
order systems of $PDE$'s, if $\mathcal{R}(\varepsilon )_{kj}^{i}(x,y)=0$ on $%
U\times V,$ then we may assign any $1$-arrow in $\varepsilon (\mathcal{U}%
_{0})$ with source at $p\in U$ and target at $q\in V$ as initial condition
and solve (39) uniquely for some $f$ on $\overline{U}\subset U$ satisfying $%
f(p)=q.$

\begin{definition}
$\varepsilon $ is flat if for any $p,q\in M$, there exist neighborhoods $%
p\in (U,x^{i}),$ $q\in (V,y^{i})$ such that $\mathcal{R}(\varepsilon
)(x,y)=0 $ on $U\times V.$
\end{definition}

It is easy to check that $\mathcal{R}(\varepsilon )(p,q)\in \Lambda
^{2}(T_{p}^{\ast })\otimes T_{q}.$ In particular, flatness has a meaning
independent of coordinates. Note that $\mathcal{R}(\varepsilon )$ is defined
on $M\times M$ and not on $M.$ Also, $\mathcal{R}(\varepsilon )(p,p)=0$
since identity diffeomorphism is the unique solution of (39). So $\mathcal{R}%
(\varepsilon )$ vanishes identically on the diagonal of $M\times M.$

Now suppose $\varepsilon $ is flat$.$ Since $\varepsilon $ preserves the
composition and inversion of $1$-arrows, the local diffeomorphisms which
integrate these $1$-arrows form a pseudogroup with the property that any
local diffeomorphism of this pseudogroup is determined on its domain by any
of its $1$-arrows. We denote this pseudugroup by $\mathcal{S}_{\varepsilon
}. $ Clearly, if $f\in \mathcal{S}_{\varepsilon },$ then all $1$-arrows of $%
f $ belong to $\varepsilon (\mathcal{U}_{0})$ by (39).

Now let $f_{1}\in \mathcal{S}_{\varepsilon }$ with domain $U_{1}$ and $%
\alpha (t),$ $0\leq t\leq 1$ a (continuous) path in $M$ with $\alpha (0)\in
U_{1}.$ Suppose there exist open sets $U_{i},1\leq i\leq k,$ which cover the
image of $\alpha $ and $f_{i}\in \mathcal{S}_{\varepsilon }$ with domain $%
U_{i}$ such that $f_{i}=f_{j}$ on $U_{i}\cap U_{j}.$ We will call this data
a continuation of $f_{1}$ along $\alpha $. The continuation of $f\in 
\mathcal{S}_{\varepsilon }$ along $\alpha $ is unique with the obvious
meaning of uniqueness. Since $f\in \mathcal{S}_{\varepsilon }$ is determined
by its $1$-arrows, continuing $f$ is the same as continuing its $1$-arrow
with source at $\alpha (0).$

\begin{definition}
Suppose $\varepsilon $ is flat$.$ If all elements of $\mathcal{S}%
_{\varepsilon }$ can be continued along paths, then $\varepsilon $ is
complete.
\end{definition}

\begin{lemma}
If $M$ is compact, then $\varepsilon $ is complete.
\end{lemma}

Proof: Let $f\in \mathcal{S}_{\varepsilon }$ with domain $U$ and $\alpha (t)$
be a path with $\alpha (0)=p\in U.$ Let $f(p)=q$ and $f^{p,q}$ denote the $1$%
-arrow of $f.$ Therefore $f^{p,q}$ has a (unique) continuation $f^{\alpha
(t),f(\alpha (t))}$ for small $t.$ Let $t_{0}\overset{def}{=}\sup \{t\in
(0,1]$ $\mid $ $f^{p,q}$ has a continuation on $(0,t)\}.$ Let $\beta (t)$
denote the curve traversed by the targets of the $1$-arrows as their source
is continued along $\alpha .$ Since $\beta $ is locally the image of $\alpha 
$ by a local diffeomorphism belonging to $\mathcal{S}_{\varepsilon }$, $%
\beta $ is continuous. Letting $t_{n}\rightarrow t_{0}^{-}$, $\beta
(t_{n_{k}})\rightarrow r\in M$ since $M$ is compact. Since $\beta $ is
continuous, we conclude $\beta (t)\rightarrow r$ as $t\rightarrow t_{0}^{-}.$
Thus we obtain a $1$-arrow $g^{\alpha (t_{0}),r}$ which belongs to $%
\varepsilon (\mathcal{U}_{0})$ since $\varepsilon (\mathcal{U}_{0})\subset 
\mathcal{U}_{1}$ is closed. Since $\mathcal{R}(\varepsilon )=0$, we can
assign $g^{\alpha (t_{0}),r}$ as initial condition and solve (39) uniquely
for some $g\in \mathcal{S}_{\varepsilon }$ such that the $1$-arrow of $g$ is 
$g^{\alpha (t_{0}),r}$ which will give a contradiction unless $t_{0}=1.$ \ $%
\square $

Now let $\widetilde{G}$ be any connected Lie group and $\Delta \subset 
\widetilde{G}$ a discrete subgroup. We call the (left) coset space $%
\widetilde{G}/\Delta $ a discrete quotient.

The geometric meaning of flatness is clarified by the following

\begin{proposition}
If $M$ is a discrete quotient, then $M$ posesses a canonical flat $%
\varepsilon $. Conversely, let $\varepsilon $ be a flat splitting$.$ If $%
\varepsilon $ is complete, then $M$ is a discrete quotient.
\end{proposition}

Proof: For the first statement, let $\pi :\widetilde{G}\rightarrow 
\widetilde{G}/\Delta =M$ be the covering map. It is well known that $M$ is
parallelizable. In fact, let $p,q\in M$. Since the action of $\Delta $ on $%
\widetilde{G}$ commutes with $\pi $, any $1$-arrow on $\widetilde{G}$
(induced by the action of some $g\in \widetilde{G}$ ) with source in $\pi
^{-1}(p)$ and target in $\pi ^{-1}(q)$ projects to the same $1$-arrow on $M$
from $p$ to $q.$ So we get a splitting $\varepsilon $ on $M$. Now $\mathcal{R%
}(\varepsilon )=0$ because the action of $\widetilde{G}$ on itself projects
locally to $M$ giving a pseudogroup on $M$ whose local diffeomorphisms
integrate these $1$-arrows on $M$ (clearly, this pseudugroup is $\mathcal{S}%
_{\varepsilon })$.

For the converse, let $\pi :\widetilde{M}\rightarrow M$ be the universal
covering map. $\mathcal{S}_{\varepsilon }$ pulls back to a pseudogroup $%
\widetilde{\mathcal{S}}_{\varepsilon }$ on $\widetilde{M}$ with the same
property: any $\widetilde{f}\in \widetilde{\mathcal{S}}_{\varepsilon }$ is
determined on its domain by any of its $1$-arrows. Let $\widetilde{f}\in $ $%
\widetilde{\mathcal{S}}_{\varepsilon }$ and let $\widetilde{\alpha }(t),$ $%
0\leq t\leq 1,$ be a path with $\widetilde{\alpha }(0)\in Dom(\widetilde{f}%
). $ Now $\widetilde{f}$ and $\widetilde{\alpha }(t)$ project to $f\in 
\mathcal{S}_{\varepsilon }$ and $\alpha (t)$ (restricting $Dom(\widetilde{f}%
) $ if necessary) and we can translate $f$ along $\alpha $ by Definition 13,
obtaining some $g\in \mathcal{S}_{\varepsilon }$ with $\alpha (1)\in Dom(g).$
Lifting this \textquotedblleft analytic continuation\textquotedblright\ to $%
\widetilde{M},$ we obtain some $\widetilde{g}\in \widetilde{\mathcal{S}}%
_{\varepsilon }$ with $\widetilde{\alpha }(1)\in Dom(\widetilde{g}).$ Since $%
\widetilde{M}$ is simply connected, the standard monodromy argument shows
that this analytic continuation is independent of the choice of the path
from $\widetilde{\alpha }(0)$ to $\widetilde{\alpha }(1).$ Therefore, any $%
\widetilde{f}\in \widetilde{\mathcal{S}}_{\varepsilon }$ extends uniquely to
a global diffeomorphism $\widetilde{f}^{e}$ on $\widetilde{M}.$ We define $%
\widetilde{G}\overset{def}{=}\{\widetilde{f}^{e}\mid \widetilde{f}\in 
\widetilde{\mathcal{S}}\}$ and check that $\widetilde{G}$ is a group.
Clearly $\widetilde{G}$ acts simply transitively on $\widetilde{M}$. In
fact, it is not difficult to show that $\widetilde{G}$ is a Lie group so
that $\widetilde{M}$ is the underlying manifold of $\widetilde{G}.$ Letting $%
\Delta =$ $($deck transformations on $\widetilde{M})\simeq \pi _{1}(M),$ we
get $\widetilde{G}/\Delta =M$ and $M$ is a discrete quotient. \ $\square $

The next proposition will add another equivalent condition to the conditions
of Proposition 8 which we denoted by $\mathbf{A.}$

\begin{proposition}
$\mathcal{R}(\varepsilon )=0$ if and only if $\mathbf{A}$ holds.
\end{proposition}

Proof: Suppose $\mathcal{R}(\varepsilon )=0$ on $M\times M$. We
differentiate the LHS of (40) with respect to $y$ and set $y=x$, which gives 
$\mathcal{R}_{2}(\varepsilon )=0$ after some computation using (11). For the
converse, we first observe that $\mathcal{R}(\varepsilon )=0$ on $M\times M$
if and only if $\mathcal{R}(\varepsilon )=0$ on some neighboorhood $U\times
U $ of $(p,p)$ for any $p\in M.$ To see this, let $p,q\in M$ and consider
the $1$-arrow $\varepsilon ^{p,q}.$ We choose \textit{any }path from $p$ to $%
q$ and take a covering $(U_{i})$ of this path such that $\mathcal{R}%
(\varepsilon )=0$ on $U_{i}\times U_{i}$. We subdivide this path into small
paths such that each small path is contained in some $U_{i}.$ Let $p_{k},$ $%
p_{k+1}$ be the initial and terminal points of the $k$'th small path. The
composition of the $1$-arrows $\varepsilon ^{p_{k},p_{k+1}}$ is $\varepsilon
^{p,q}$ and each $1$-arrow $\varepsilon ^{p_{k},p_{k+1}}$ integrates by
assumption to a local diffeomorphism which is a solution of (39). Now the
composition of these local diffeomorphisms integrates $\varepsilon ^{p,q}$
and therefore is a solution of (39) which shows $\mathcal{R}(\varepsilon
)(p,q)=0$, proving the claim.

Now suppose $\mathcal{R}_{2}(\varepsilon )=0$ and let $(U,x^{i})$ be any
coordinate patch. We choose some $p\in U,$ a basis of the tangent space at $%
p $ and extend this basis to $\varepsilon $-invariant vector fields $\xi
_{(i)} $ on $U,$ $1\leq i\leq n.$ Now $[\xi _{(i)}(x),\xi
_{(j)}(x)]=c_{ij}^{a}(x)\xi _{(a)}(x)$ for some functions $c_{ij}^{a}(x)$ on 
$U.$ Since $\mathcal{R}_{2}(\varepsilon )=0$, the functions $c_{ij}^{a}(x)$
are constant by Proposition 8. By Lie's Third Fundamental Theorem $(LTFT)$,
there exists a local Lie group acting on $U$ (shrinking $U$ if necessary)
with infinitesimal generators $\xi _{(i)}.$ Now the proof of $LTFT$ (see
[25], pg. 398-415) shows that the $1$-arrows of the local diffeomorphisms
induced by this local Lie group action coincide with the $1$-arrows of $%
\varepsilon (\mathcal{U}_{0}).$ Thus $\mathcal{R}(\varepsilon )=0$ on $%
U\times U$ and therefore on $M\times M$ since $U$ is arbitrary$.$ \ $\square 
$

In view of Propositions 8, 16, a local Lie group is a differentiable
manifold $M$ \textit{together }with a flat splitting $\varepsilon .$ We do
not know whether a parallelizable manifold $M$ can be a local Lie group in
different ways, that is, whether it can admit two flat splittings $%
\varepsilon ,\overline{\varepsilon }$ with nonisomorphic Lie algebras $%
\mathfrak{X}_{\varepsilon }(M$ $),$ $\mathfrak{X}_{\overline{\varepsilon }%
}(M).$ Henceforth, the statement "$M$ is a local Lie group" will always
refer to $M$ together with its flat splitting which we fix once and for all.

Proposition 8 shows that the condition $\mathcal{R}_{2}(\varepsilon )=0$ is
equivalent to the hypothesis of $LTFT$ and the proof of Proposition 16 shows
that $\mathcal{R}(\varepsilon )=0$ is equivalent to its conclusion.
Therefore, the program outlined in the Introduction, if viewed locally, is
nothing but a natural generalization of of $LTFT$ from simply transitive
actions (where the geometric order $m=0$) to transitive and effective
actions (where $m$ is arbitrary).

As we see from the proof of Proposition 16, $\widetilde{\mathcal{S}}%
_{\varepsilon }$ globalizes to the Lie group $\widetilde{G}$ with underlying
manifold $\widetilde{M}$. However, $\mathcal{S}_{\varepsilon }$ may
globalize to a Lie group $G$ already on $M$, that is, $\varepsilon $ may
have "trivial monodromy" already on $M$ (on the other hand, recall that $%
\varepsilon ,$ as the trivialization of the principal bundle $\mathcal{U}%
_{1}^{(p)}$ $\rightarrow M,$ has necessarily trivial monodromy).

\begin{definition}
A local Lie group $M$ is globalizable if any $f\in \mathcal{S}_{\varepsilon
} $ extends (necessarily uniquely) to a global diffeomorphism of $M.$
\end{definition}

So $M$ is globalizable if and only if there are no obstructions to the
existence of \textit{global }solutions of (39) in which case the \textit{%
global }diffeomorphisms of $\mathcal{S}_{\varepsilon }$ act simply
transitively on $M$ and have the structure of a Lie group $G.$ In this case,
we simply say that $M$ is a Lie group. Observe that the statement
\textquotedblleft $M$ is a Lie group\textquotedblright\ refers to two
different objects: the pair $(G,M)$ where $G$ is a transformation group of $%
M $ which acts simply transitively and the abstract Lie group $G.$ We can
identify $G$ and $M$ by choosing some point $e\in M$ and map $g\in G$ to its 
$0$-arrow from $e$ to $g(e).$ Clearly, this map $G\rightarrow M$ is a
bijection, in fact, a diffeomorphism, but there is no canonical
identification. We will continue to make such identification as we did
already in the proof of Proposition 15.

Definition 17 of globalizability coincides with the one given in [20] if $M$
is complete.

Using the notation and setting of Proposition 15, the next proposition
characterizes globalizability of $M$ in terms of the discrete subgroup $%
\Delta \subset \widetilde{G}.$

\begin{proposition}
$M$ is a Lie group if and only if $\Delta \subset \widetilde{G}$ is a normal
(hence central) subgroup.
\end{proposition}

Proof: If $\Delta \subset \widetilde{G}$ is a normal subgroup, then $\pi :%
\widetilde{G}\longrightarrow \widetilde{G}/\Delta =M$ is a homomorphism of
Lie groups and surely $M$ is globalizable.

For the converse, suppose that $M$ is globalizable. Let $\pi (e)=x,$ $g\in 
\widetilde{G},$ $\pi (g)=y.$ Clearly $\pi (d(0))=x$ for all $d(0)\in \Delta
. $ Now the action of $g$ on $\widetilde{G}$ defines a local diffeomorphism
which maps a neighborhood of the identity $e$ to a neighborhood of $g.$ This
local diffeomorphism projects to some $f(1)\in \mathcal{S}_{\varepsilon }$
with $f(1)(x)=y.$ We choose a path $A(s),$ $0\leq s\leq 1$ with $A(0)=e,$ $%
A(1)=d(0)$ and consider the loop $\pi \circ A=\overline{A}$ at $x.$ We now
continue $f(1)$ along the loop $\overline{A}$ turning the source $x$ to its
initial value. Thus the target $f(1)(x)=y$ traces a curve starting from $y.$
Since $f(1)$ globalizes by assumption, $y$ also turns to its initial value
and this curve is also a loop at $y.$ We lift this loop to a curve $B(s),$ $%
0\leq s\leq 1,$ which starts from $B(0)=g.$ Since $\pi (B(1))=\pi (B(0))=y$,
there exists a unique $d(1)\in \Delta $ with $B(1)=d(1)g.$ We now lift the
continuation of $f(1)\in \mathcal{S}_{\varepsilon }$ to $\widetilde{G}$ as
in the proof of Proposition 15$.$ Since the action of $g$ is a global
diffeomorphism and $g$ remains the same during this continuation, we
conclude $gA(s)=B(s)$ for all $s.$ Therefore $gd(0)=gA(1)=gB(1)=d(1)g.$ We
claim that $d(0)=d(1).$ This will finish the proof since $g$ and $d(0)$ are
both arbitrary in $gd(0)=d(0)g.$

The prove the claim, we choose a path $g(t)$ with $g(0)=e,$ $g(1)=g.$
Replacing $g=g(1),$ $f(1)\in \mathcal{S}_{\varepsilon }$ with $g(t),$ $f(\pi
(g(t))\in $ $\mathcal{S}_{\varepsilon }$ and repeating the above argument,
we obtain a loop at $\pi (g(t)$ which lifts to a curve $H(s,t),$ $0\leq
s\leq 1$ with $H(0,t)=g(t).$ Clearly, $H(s,0)=A(s)$ and $H(s,1)=B(s).$ By
construction, $H(s,t)$ is continuous on $[0,1]\times \lbrack 0,1]$ and
therefore $H(1,t)$ is a curve from $d(0)$ to $d(1)g.$ Since $\pi
(H(0,t))=\pi (H(1,t)),$ there exists a unique $d(t)\in \Delta $ satisfying $%
d(t)g(t)=H(1,t).$ We define $F(t)\overset{def}{=}$ $H(1,t)g(t)^{-1}\in
\Delta $ for all $t.$ Since $\Delta $ is discrete and $F(t)$ is continuous, $%
F(t)$ is constant which proves the claim. \ $\square $

Recalling Definition 10, a Lie group is clearly a complete local Lie group.
The proof of the next corollary is immediate from the proof of Proposition
15.

\begin{corollary}
A simply connected and complete local Lie group is a Lie group. \ \ 
\end{corollary}

Henceforth in this section, we assume that $M$ is a complete local Lie group
and therefore a discrete quotient in view of Proposition 16.

Now, since the $1$-form $w$ is closed by Corollary 9, it is locally exact
and our aim is to find an explicit local primitive of $w.$ For this purpose,
we fix some $e\in M$ and a simply connected coordinate patch $e\in
(U,x^{i}). $ For $a,b\in U$, we define $ab\in M$ as follows: since $\mathcal{%
R}(\varepsilon )=0,$ the $1$-arrows $\varepsilon ^{e,a}$ and $\varepsilon
^{e,b}$ integrate to local diffeomorphisms $f,g\in \mathcal{S}_{\varepsilon
} $ with $f(e)=a$ and $g(e)=b$ and $f,g$ are both defined on $U$ since $U$
is simply connected and $\varepsilon $ is complete. We define

\begin{equation}
ab\overset{def}{=}(f\circ g)(e)=f(b)\in M
\end{equation}
\ 

In (41) we adhere to the standard convention of composing from right to
left. Note that (41) depends on our choice of the base point $e$ and $ab$
may not be in $U$ unless $a,$ $b$ are sufficiently close to $e$ $.$ At this
point, it is worthwhile to observe that the classical concept of local Lie
group is modelled on this local structure, that is, the set $U$ together
with the local multiplication (41), whereas a local Lie group $M$ according
to Definition 10 is essentially a global concept. In particular, it may be
misleading to imagine a local Lie group $M$ as an open neighborhood of
identity in some global Lie group $G.$ In fact, the reader may have observed
that it is possible to define a Lie group as a globalizable local Lie
group.\ \ 

The local multiplication (41) determines the local left and right
translations $L_{a},$ $R_{a}$ defined by $L_{a}(x)=ax$ and $R_{a}(x)=xa$
where $a,x\in U.$ Thus we may consider the $1$-arrows of $L_{a}$ and $R_{a}$
with source at $e$ and target at $a$, that is $j_{1}(L_{a})^{e,a}$ and $%
j_{1}(R_{a})^{e,a}$ which we want to compute now in coordinates.

Let $f\in \mathcal{S}_{\varepsilon }$, $f(e)=a$ and $f$ be defined on $U$.
Also, for any $x\in U$, let $g_{x}\in \mathcal{S}_{\varepsilon }$ such that
the diffeomorphism $y\rightarrow g_{x}(y)$ is defined on $U$ and is the
unique solution of (39) in the variable $y$ satisfying the initial condition 
$g_{x}(e)=x.$ We also have $g_{e}(y)=y$ since the only solution of (39)
which fixes $e$ is identity on $U$. Now (41) gives $L_{a}(x)=ax=(f\circ
g_{x})(e)=f(x)$ and therefore

\begin{equation}
\left[ j_{1}(L_{a})^{e,a}\right] _{j}^{i}=\left[ \frac{\partial f^{i}(x)}{%
\partial x^{j}}\right] _{x=e}=\varepsilon _{j}^{i}(e,a)
\end{equation}

We see from (42) that the assumption $\mathbf{A}$ is not needed to define $%
j_{1}(L_{a})^{e,a}$ and the "local left $1$-arrows" are actually global as
they coincide with the $1$-arrows of $\varepsilon $.

Similarly

\begin{equation}
\left[ j_{1}(R_{a})^{e,a}\right] _{j}^{i}=\left[ \frac{\partial ((g_{x}\circ
f)(e))^{i}}{\partial x^{j}}\right] _{x=e}=\left[ \frac{\partial g_{x}^{i}(a)%
}{\partial x^{j}}\right] _{x=e}
\end{equation}

We see from (43) that the definition of $j_{1}(R_{a})^{e,a}$ depends
essentially on $\mathbf{A}$ and "right $1$-arrows" are local and defined
only on $U.$ The following points are worth mentioning: if $a\in M$, $b\in U$%
, then $ab$ is still defined by (41), whereas to define $ba$ we must
continue the domain of $f$ so that it contains the point $b,$ but
continuations along different paths may give different values for $f(b)$ and
it may be impossible to define $ba$ uniquely. Equivalently, we have the
multiplication map $M\times U\rightarrow M$ but not necessarily $U\times
M\rightarrow M,$ in contrast to the classical definition of a local Lie
group. It is due to this fact that left $1$-arrows are global whereas right $%
1$-arrows are not. Also, left and right are reversed if we change our
convention of composition from right to left in (41), which shows that left
and right are nonconcepts for the local Lie group $M$ (even if it is
globalizable!) whereas they are honest concepts for the Lie group $G.$

We now define the function $Ad_{e}:U\rightarrow $ $\mathcal{U}_{1}^{e,e}$ by

\begin{equation}
Ad_{e}(a)\overset{def}{=}\left[ j_{1}(L_{a})^{e,a}\right] ^{-1}\circ
j_{1}(R_{a})^{e,a}
\end{equation}

We call $Ad_{e}$ the local adjoint map based at $e.$ Using (42) and (43) we
get

\begin{equation}
\left[ Ad_{e}(a)\right] _{j}^{i}=\varepsilon _{b}^{i}(a,e)\left[ \frac{%
\partial g_{x}^{b}(a)}{\partial x^{j}}\right] _{x=e}
\end{equation}

Since $Ad_{e}(e)=id^{e,e}$, $\det (Ad_{e}(x))$ is positive at $x=e$ and
therefore positive on $U$ since it is nonzero and continuous.

\begin{lemma}
Let $A(x)=\left[ \alpha _{j}^{i}(x)\right] $ be an $n\times n$ matrix whose
entries are smooth functions on an open set $U\subset \mathbb{R}^{n}.$
Suppose that $A(x)$ is invertible on $U$ with inverse $B(x)=\left[ \beta
_{j}^{i}(x)\right] .$ Then

\begin{equation}
(\det A)^{-1}\frac{\partial \det (A)}{\partial x^{j}}=\frac{\partial \alpha
_{b}^{c}}{\partial x^{j}}\beta _{c}^{b}
\end{equation}
\end{lemma}

\QTP{Body Math}
Proof: (46) follows easily from the cofactor expansion of $\det (A)$ (see
[25], pg. 8, formula 7.2). \ $\square $

\QTP{Body Math}
Note that the $RHS$ of (46) is $Tr\left( \frac{\partial A}{\partial x^{j}}%
A^{-1}\right) $. If $\det (A)$ is positive on $U,$ then the $LHS$ of (46) is
equal to

\QTP{Body Math}
\begin{equation*}
\frac{\partial \log (\det (A))}{\partial x^{j}}
\end{equation*}

\begin{proposition}
$-d(\log (\det Ad_{e}))=w$ on $U.$
\end{proposition}

Proof: We first show that $-d(\log (\det Ad_{e}))$ and $w$ have the same
values at $x=e.$ Using (45), Lemma 20, (11) and (36) we compute

\begin{eqnarray}
-\left[ d(\log (\det Ad_{e}))(e)\right] _{i} &=&-\left[ \frac{\partial }{%
\partial y^{i}}\log \det \left( \varepsilon _{b}^{j}(y,e)\frac{\partial
g_{x}^{b}(y)}{\partial x^{k}}\right) \right] _{x=e,y=e}  \notag \\
&=&-\left[ \frac{\partial }{\partial y^{i}}\left( \varepsilon _{b}^{c}(y,e)%
\frac{\partial g_{x}^{b}(y)}{\partial x^{d}}\right) \right]
_{x=e,y=e}(\delta _{c}^{d})  \notag \\
&=&-\left[ \frac{\partial }{\partial y^{i}}\left( \varepsilon _{b}^{c}(y,e)%
\frac{\partial g_{x}^{b}(y)}{\partial x^{c}}\right) \right] _{x=e,y=e} 
\notag \\
&=&-\left[ \frac{\partial \varepsilon _{b}^{c}(y,e)}{\partial y^{i}}\frac{%
\partial g_{x}^{b}(y)}{\partial x^{c}}+\varepsilon _{b}^{c}(y,e)\frac{%
\partial ^{2}g_{x}^{b}(y)}{\partial y^{i}\partial x^{c}}\right] _{x=e,y=e} 
\notag \\
&=&-\left[ \frac{\partial \varepsilon _{b}^{c}(y,e)}{\partial y^{i}}\right]
_{y=e}(\delta _{c}^{b})-(\delta _{b}^{c})\frac{\partial }{\partial x^{c}}%
\left[ \frac{\partial g_{x}^{b}(y)}{\partial y^{i}}\right] _{x=e,y=e}  \notag
\\
&=&\Gamma _{ib}^{b}(e)-\frac{\partial }{\partial x^{b}}\left[ \frac{\partial
g_{x}^{b}(y)}{\partial y^{i}}\right] _{x=e,y=e}  \notag \\
&=&\Gamma _{ib}^{b}(e)-\left[ \frac{\partial \varepsilon _{i}^{b}(e,x)}{%
\partial x^{b}}\right] _{x=e}  \notag \\
&=&\Gamma _{ib}^{b}(e)-\Gamma _{bi}^{b}(e)  \notag \\
&=&\left[ w(e)\right] _{i}
\end{eqnarray}%
which proves the claim. Now let $\overline{e}\in U$ another base point. It
is easy to check that

\begin{equation}
\log (\det Ad_{e}(x))=\log (\det Ad_{\overline{e}}(x))+\log (\det Ad_{e}(%
\overline{e}))
\end{equation}

Differentiating (48) at $x=\overline{e}$ and using (47), we get

\begin{equation}
-d(\log (\det Ad_{e})(\overline{e})=w(\overline{e})
\end{equation}%
which finishes the proof since $\overline{e}$ is arbitrary. \ $\square $

Note that (44) defines $Ad$ by the formula $x\rightarrow g^{-1}xg$ rather
than the standard one $x\rightarrow gxg^{-1}$ which accounts for the minus
sign in Proposition 21. The reason for our choice is that it was easier to
invert the local formula for $j_{1}(L_{a})^{e,a}$.

Now let $M$ be a Lie group, that is, $M$ be globalizable. We fix some $e\in
M $ and define the function $Ad:M\rightarrow \mathcal{U}_{1}^{e,e}$ by (44).
Observe that $Ad$ is now defined on $M$ since right $1$-arrows are globally
defined since $M$ is globalizable. Thus we have

\begin{proposition}
$w$ is exact on a Lie group $M$ with primitive $-d(\log (\det Ad_{e})).$
\end{proposition}

We refer the reader to [15], Propositions 6.29, 6.34 for an intriguing
relation of the primitive $-\log (\det Ad_{e})$ to the Maslov index. See
also [9] for another interpretation of $Ad_{e}$ in $-\log (\det Ad_{e})$ as
a Nijenhuis tensor.

\begin{proposition}
Let $f:M\rightarrow N$ be a smooth covering map. If $N$ is a local Lie
group, then $f$ determines a local Lie group structure on $M$ and $f^{\ast
}w_{N}=w_{M}.$
\end{proposition}

Proof: Straightforward using definitions. \ $\square $

\begin{proposition}
Let $\widetilde{G}$ be a connected Lie group, $\Delta \subset $ $\widetilde{G%
}$ a discrete subgroup and $\pi :$ $\widetilde{G}\rightarrow \widetilde{G}%
/\Delta =M$ the covering map. Then $w_{M}$ is exact if and only if $\log
(\det (Ad_{\widetilde{G}}))(ax)$ $=x,$ $a\in $ $\Delta ,$ $x\in \widetilde{G}%
,$ that is, the function $\log (\det (Ad_{\widetilde{G}}))$ is automorphic
with respect to $\Delta .$
\end{proposition}

Proof: $\log (\det (Ad_{\widetilde{G}}))$ is constant on the orbits of $%
\Delta $ if and only if there exists a function $h$ on $M$ $=\widetilde{G}%
/\Delta $ with $h\circ \pi =\log (\det (Ad_{\widetilde{G}})).$ We choose an
open set $U\subset H$ such that $\pi :U\rightarrow V\subset M$ is a
diffeomorphism. Now $w_{\widetilde{G}}=d\log (\det (Ad_{\widetilde{G}%
}))=d(h\circ \pi )=\pi ^{\ast }(dh)$ on $U.$ Therefore $(\pi ^{\ast
})^{-1}(w_{\widetilde{G}})=w_{M}=dh$ on $V$ by Lemma 23$.$ Since $V$ is
arbitrary, we conclude $dh=w_{M}$ on $M$ and the conclusion follows. \ $%
\square $

Recall that a connected Lie group is called unimodular if $\det (Ad)=1.$ $G$
is unimodular if and only if it admits a nonzero invariant $n$-form.

\begin{corollary}
If $M$ is the discrete quotient of a connected and unimodular Lie group $%
\widetilde{G}$, then $w_{M}=0.$
\end{corollary}

Proof: By Lemma 23 and Proposition 22, we have $\pi ^{\ast }w_{M}=w_{%
\widetilde{G}}=d\log (\det Ad)=0$. Therefore $w_{M}=0$ since $\pi $ is a
local diffeomorphism. \ $\square $

Many important Lie groups are unimodular: abelian, nilpotent, compact,
semisimple and reductive. Now let $B\subset SL(2,\mathbb{R)}$ the Borel
subgroup of upper triangular matrices. An easy computation shows $\det (ad%
\left[ 
\begin{array}{cc}
a & b \\ 
0 & c%
\end{array}%
\right] )=\frac{c}{a}.$ Let $\Delta \subset B$ be the subgroup of matrices
of the form $\left[ 
\begin{array}{cc}
2^{n} & 0 \\ 
0 & 2^{-n}%
\end{array}%
\right] $ where $n$ is an integer. Clearly $\Delta $ is a discrete subgroup
and $\det (Ad)$ is not constant on $\Delta $ so that $[w_{\frac{B}{\Delta }%
}]\neq 0.$

In Definition 41 and therefore in Proposition 21 we assumed completeness of $%
\varepsilon $. However this assumption can be dropped since we can choose $U$
sufficiently small in Definition 41 so that \textit{all }$1$-arrows with
source and target in $U$ integrate to local diffeomorphisms defined on $U.$
Therefore $[w]$ is in force also if $M$ is not complete. Many explicit and
nontrivial examples of incomplete local Lie groups are constructed in [20]
which fail to be globalizable. The example on pg. 49 in [20] is even simply
connected, but such pathalogy can not occur if $M$ is complete in view of
Corollary 19. We will postpone the study of these examples to some future
work.

We will conclude this section with three remarks.

1) If $\varepsilon $ is complete, $\mathbf{A}$ forces $M$ to be a discrete
quotient by Proposition 15 and hence analytic. Therefore $\mathbf{A}$ is a
strong condition. On the other hand, we only need $Tr\mathcal{R}%
_{2}(\varepsilon ))=0$ to prove $dw=0$ and not $\mathbf{A.}$ However, we are
unable to attach a seperate geometric meaning to the condition $Tr\mathcal{R}%
_{2}(\varepsilon )=0.$

2) We fix some $e\in M$ and a coordinate system around $e$ once and for all
and write $\varepsilon _{j}^{i}(x,y)=\varepsilon _{a}^{i}(e,y)\varepsilon
_{j}^{a}(x,e).$ Using the notation of [23], [24], we define a geometric
object $w$ (not to be confused with the $1$-form $w)$ with components $%
w_{j}^{i}(x)$ over $(U,x^{i})$ defined by $w_{j}^{i}(x)\overset{def}{=}%
\varepsilon _{j}^{i}(x,e)$ so that we have

\begin{equation}
\varepsilon _{j}^{i}(x,y)=\overline{w}_{a}^{i}(y)w_{j}^{a}(x),\text{ \ \ \ \
\ }\overline{w}_{a}^{i}(y)w_{j}^{a}(y)=\delta _{j}^{i}
\end{equation}

Replacing derivatives with jet sections in the transformation rule $%
w_{j}^{i}(x)=w_{a}^{i}(y)\frac{\partial y^{a}}{\partial x^{i}}$ gives

\begin{equation}
w_{j}^{i}(x)=w_{a}^{i}(y)f_{i}^{a}(x)
\end{equation}%
which describes the association of the object $w$ with the groupoid $%
\mathcal{U}_{1}.$ Now (51) is taken as the basis of the Maurer-Cartan form
in [23], [24] (see [23], pg. 212-221, 316-317, [24], pg. 27-32, 246-249).
Substituting (51) into (40), a straightforward computation shows that the
condition $\mathcal{R}(x,y)=0$ is equivalent to the second formula on pg. 28
in [24] which is

\begin{equation}
\overline{w}_{j}^{a}(y)\overline{w}_{k}^{b}(y)\left( \frac{\partial
w_{b}^{i}(y)}{\partial y^{a}}-\frac{\partial w_{a}^{i}(y)}{\partial y^{b}}%
\right) =\overline{w}_{j}^{a}(x)\overline{w}_{k}^{b}(x)\left( \frac{\partial
w_{b}^{i}(x)}{\partial x^{a}}-\frac{\partial w_{a}^{i}(x)}{\partial x^{b}}%
\right)
\end{equation}

Thus both sides of (52) must be equal to the same constants $c_{jk}^{i}$. It
is shown in [23], [24] that (51) is formally integrable if and only if (52)
holds. Thus we see that (51) is a consequence of the existence of the
splitting $\varepsilon $ and and the formal integrability of (51), a concept
which plays a fundamental role in the works [23], [24], is equivalent to the
condition $\mathcal{R}(\varepsilon )=0$.

Note that (16), (18) and (50) give%
\begin{equation}
\Gamma _{kj}^{i}(y)dy^{k}=\frac{\partial \overline{w}_{a}^{i}(y)}{\partial
y^{k}}w_{j}^{a}(y)dy^{k}
\end{equation}

The $1$-form on the LHS of (53) has only local meaning. On the other hand, $%
\overline{w}_{j}^{i}(y)=\varepsilon _{j}^{i}(y,e)$ uniquely determines some $%
g\in \mathcal{S}_{\varepsilon }$ in view of (41) and therefore we may write
the RHS of (53) as $(dg)g^{-1}.$

3) Let $M$ be parallelizable with the splitting $\varepsilon $. We fix $e\in
(U,x^{i})$, a tangent vecur $X_{e}$ at $e,$ and define the $ODE$ for the
unknown curve $x(t)$, $0\leq t\leq \epsilon ,$ $x(0)=e,$ by $(\varepsilon
^{e,x})_{\ast }(X_{e})=$ $\overset{\cdot }{x}(t)$ for $x\in U.$ In
coordinates,

\begin{equation}
\frac{dx^{i}}{dt}=\varepsilon _{a}^{i}(e,x)X^{a},\text{ \ }1\leq i\leq n
\end{equation}

We call local solutions of (54) left $1$-parameter curves. Differentiating
(54) at $x=e$ and using (11), we get

\begin{equation}
\frac{d^{2}x^{i}}{dt^{2}}(e)-\Gamma _{ab}^{i}(e)\frac{dx^{b}}{dt}(e)\frac{%
dx^{a}}{dt}(e)=0
\end{equation}%
which are the equations for geodesics (note again the minus sign in (55)
which arises from our convention of differentiating with respect to the
second argument in (11)). Thus we can define left geodesic completeness of $%
M.$ If $M$ is a local Lie group, then we can define locally also right $1$%
-parameter curves and therefore $1$-parameter subgroups. It is an
interesting problem to study the relations between these well known objects
(which is done to some extent in [14], pg. 26-47, 97-131) and the relation
between the concepts of completeness according to Definition 13 and geodesic
completeness.

\section{Odd degree forms}

In this section we will define the higher order analogs of the closed $1$%
-forms $\omega $ and $w.$

Turning back to (29), let $X_{1},...,X_{k}$ be sections of $J_{1}T.$ We
define the $k$-form $\omega _{k}$ of $J_{1}T$ by

\begin{equation}
\omega _{k}(X_{1},...,X_{k})\overset{def}{=}\frac{1}{k}\sum_{\sigma \in
S_{k}}sgn(\sigma )Tr(X_{\sigma (1)}\circ X_{\sigma (2)}\circ ...\circ
X_{\sigma (k)})
\end{equation}

Clearly, $\omega _{1}=\omega $.

\begin{proposition}
$\omega _{2m}=0$, $m\geq 1$
\end{proposition}

\bigskip Proof: Consider the set $K$ of $k$-tuples $(\sigma (1),\sigma
(2),...,\sigma (k))$ where $\sigma $ is a permutation. We define an
equivalence relation $\sim $ on $K:$ $(\sigma (1),\sigma (2),...,\sigma
(k))\sim (\tau (1),\tau (2),...,\tau (k))$ if $(\sigma (1),\sigma
(2),...,\sigma (k))=(\tau (i),\tau (i+1),...\tau (k),\tau (1),\tau (2),...,$ 
$\tau (i-1))$ for some $i,$ that is, two permutations are equivalent if they
differ by a cyclic permutation. This equivalence relation is imposed by $%
Tr(ab)=Tr(ba).$ Now consider the formal sum

\begin{equation}
\sum_{\sigma \in S_{k}}sgn(\sigma )((\sigma (1),\sigma (2),...,\sigma (k))
\end{equation}

The proof reduces to the following combinatorial statement: If $k$ is even
and we identify the equivalent $k$-tuples in (51), then (51) vanishes. To
prove this, we choose two odd integers $a,b$ \ with $k=a+b.$ For some fixed $%
\sigma ,$ we have $sgn(\sigma )(\sigma (1),...,\sigma (a),\sigma
(a+1),...\sigma (a+b))=(-1)^{ab}sgn(\sigma )(\sigma (a+1),...\sigma
(a+b),\sigma (1),...,\sigma (a))$ and all terms in (51) cancel in pairs. \ $%
\square $

Using $1$-arrows of $\varepsilon $ we can define $\varepsilon $-invariant $k$%
-forms on $M.$ Following our convention for composition in Section 8, we
will call these forms left invariant. Exterior derivative of a left
invariant form need not be left invariant, but this is so if $\mathcal{R}%
(\varepsilon )=0$, that is, if $M$ is a local Lie group. Thus we get a
complex whose cohomology $H^{\ast }(\mathfrak{X}_{\varepsilon }(M))$
coincides with the Lie algebra cohomology of $\mathfrak{X}_{\varepsilon
}(M). $ The left invariant forms can be localized at any point $p\in M$ and
we obtain a complex at $p$ with cohomology $H_{p}^{\ast }(\mathfrak{X}%
_{\varepsilon }(M)).$ Since the evaluation map $e_{p}$ in (38) is an
isomorphism of Lie algebras, it induces an isomorphism $e_{p}^{\ast
}:H_{p}^{\ast }(\mathfrak{X}_{\varepsilon }(M))\rightarrow H^{\ast }(%
\mathfrak{X}_{\varepsilon }(M)).$ If $q\in M$ is another point, there exists
a unique $f\in \mathcal{S}_{\varepsilon }$ with $f(p)=q.$ Since the local
diffeomorphism $f$ commutes with the exterior derivative, it induces an
isomorphism $f^{\ast }:$ $H_{q}^{\ast }(\mathfrak{X}_{\varepsilon
}(M))\rightarrow H_{p}^{\ast }(\mathfrak{X}_{\varepsilon }(M))$ and we
obtain the commutative diagram

\begin{equation}
\begin{array}{ccccc}
&  & H^{\ast }(\mathfrak{X}_{\varepsilon }(M)) &  &  \\ 
& e_{q}^{\ast }\nearrow &  & \nwarrow e_{p}^{\ast } &  \\ 
H_{q}^{\ast }(\mathfrak{X}_{\varepsilon }(M)) &  & \overset{f^{\ast }}{%
\longrightarrow } &  & H_{p}^{\ast }(\mathfrak{X}_{\varepsilon }(M))%
\end{array}%
\end{equation}

Now, we pull back the forms $\omega _{2k+1}$ defined by (50) to $w_{2k+1}$
on the local Lie group $M.$ Let $X_{1},X_{2},....X_{k}$ be left invariant
vector fields on $M.$ Proposition 11, the derivation of (36) and (50) give
the formula

\begin{equation}
w_{k}(X_{1},...,X_{k})=\frac{1}{k}\sum_{\sigma \in S_{k}}sgn(\sigma
)tr\left( ad(X_{\sigma (1)})\circ ...\circ ad(X_{\sigma (k)}\right)
\end{equation}

As we noted at the end of Section 6, $w_{1}$ is left invariant and (53)
shows that $w_{2k+1}$ is left invariant for all $k\geq 0.$ Therefore the
exterior derivative of $w_{2k+1}$ can be computed algebraically.

\begin{proposition}
$dw_{2k+1}=0,$ $k\geq 0$
\end{proposition}

Proof: If $h_{1}\circ h_{2}\circ ....\circ h_{k}$ is a composition of some
linear maps and $\sigma \in S_{k},$ we define the composition

\begin{equation}
\sigma \ast (h_{1}\circ h_{2}\circ ...\circ h_{k})\overset{def}{=}h_{\sigma
(1)}\circ h_{\sigma (2)}\circ ...\circ h_{\sigma (k)}
\end{equation}%
and extend $\ast $ linearly over sums of compositions. Usind the notation
(54), we now have (omitting the constant factors), $\
dw_{k}(X_{1,}...,X_{k+1})$

\begin{eqnarray}
&=&\sum_{i\leq j+1}(-1)^{i+j}w_{k}([X_{i},X_{j}],...,\widehat{X}_{i},...,%
\widehat{X}_{j},...,X_{k+1})  \notag \\
&=&\sum_{i\leq j+1}(-1)^{i+j}\sum_{\sigma \in S_{k}}sgn(\sigma )tr\left[
\sigma \ast (ad[X_{i},X_{j}]\circ .....\circ adX_{k+1})\right]  \notag \\
&=&\sum_{i\leq j+1}(-1)^{i+j}\sum_{\sigma \in S_{k}}sgn(\sigma )tr\left[
\sigma \ast ([adX_{i},adX_{j}]\circ .....\circ adX_{k+1})\right] \\
&=&\sum_{i\leq j+1}(-1)^{i+j}\sum_{\sigma \in S_{k}}sgn(\sigma )tr\left[
\sigma \ast ((adX_{i}\circ adX_{j})\circ .....\circ adX_{k+1})\right]  \notag
\\
&&-\sum_{i\leq j+1}(-1)^{i+j}\sum_{\sigma \in S_{k}}sgn(\sigma )tr\left[
\sigma \ast ((adX_{j}\circ adX_{i})\circ .....\circ adX_{k+1})\right]  \notag
\end{eqnarray}%
Now each term in the last sum \i n (55) contains a composition of $k$ linear
maps. The value of trace does not change if we apply some $k$-cycle to its
arguments. If $k$ is odd and $\sigma \in S_{k}$ is a $k$-cycle, then $%
sgn(\sigma )=1.$ Using these facts, a straightforward computation shows that
all terms in the last sum in (55) cancel in pairs. We will omit the
details.\ $\square $

Observe that (53) and Proposition 27 define odd order characteristic classes
in the Lie algebra cohomology of \textit{any abstract Lie algebra $\mathfrak{%
g}$. }We will denote these classes by $[w_{2k+1}]_{\mathfrak{g}}\in H^{2k+1}(%
\mathfrak{g).}$ In fact, one can replace $ad$ with any representation $\rho :%
\mathfrak{g}\rightarrow V$ and still define these classes which however will
depend on the representation as in [6]. We will not pursue this more general
approach here.

If $\Delta \subset G$ a discrete subgroup, then $[w_{2k+1}]_{dR}\in
H^{2k+1}(M,\mathbb{R})$ is defined in view of Proposition 27. We also have $%
[w_{2k+1}]_{\mathfrak{X}(\varepsilon )}\in H^{2k+1}(\mathfrak{X}%
_{\varepsilon }(M)).$ It is worthwhile to observe that these cohomology
classes are intrinsic characteristic classes of the local Lie group $M$ and
are defined regardless of any foliation. If $[w_{2k+1}]_{\mathfrak{X}%
(\varepsilon )}=0,$ then clearly $[w_{2k+1}]_{dR}=0$ but the converse is
false as shown by our noncompact example in Section 8: $[w_{1}]_{\mathfrak{X}%
(\varepsilon )}=0$ if and only if $w_{1}=0$ whereas we may have $%
[w_{1}]_{dR}=0$ with nonconstant primitive $-\log (\det Ad_{e})$ and
therefore $w_{1}\neq 0.$

Now let $G=SL(2,\mathbb{R)}$, $\Delta \subset G$ a cocompact discrete
subgroup and $\mathcal{F}$ the codimension one foliation on $M$ defined by
the left invariant vector fields $X=\left[ 
\begin{array}{cc}
0 & 1 \\ 
0 & 0%
\end{array}%
\right] ,$ $H=\left[ 
\begin{array}{cc}
1 & 0 \\ 
0 & -1%
\end{array}%
\right] .$ Let $GV(\mathcal{F})\in $ $H^{3}(M,\mathbb{R})$ denote the
Godbillon-Vey class of this foliation (see [4], pg. 62-64 for details). To
distinguish between $w_{2k+1}$ on the local Lie groups $G$ and $M=G/\Delta ,$
we will use the notations $w_{2k+1}^{G}$ and $w_{2k+1}^{M}.$

\begin{proposition}
$[w_{3}^{M}]_{dR}=GV(\mathcal{F})$
\end{proposition}

Proof: It suffices to check that $w_{3}^{G}$ is a nonzero multiple of an
invariant volume form on $SL(2,\mathbb{R)}.$ Let $X,$ $H,$ $Y$ $\ $be the
basis of left invariant vector fields on $SL(2,\mathbb{R)}$ where $Y=\left[ 
\begin{array}{cc}
0 & 0 \\ 
1 & 0%
\end{array}%
\right] $ and $X,$ $H$ are as above. We define the invariant volume form $%
\nu $ on $SL(2,\mathbb{R)}$ by $\nu =X^{\ast }\Lambda H^{\ast }\Lambda
Y^{\ast }$ where $X^{\ast },$ $H^{\ast },$ $Y^{\ast }$ are the dual $1$%
-forms. On the other hand, (53) gives

\begin{eqnarray}
w_{3}^{G}(X,H,Y) &=&\frac{1}{3}\{Tr(adX\circ adH\circ adY)+Tr(adY\circ
adX\circ adH)  \notag \\
&&+Tr(adH\circ adY\circ adX)-Tr(adX\circ AdY\circ adH)  \notag \\
&&-Tr(adH\circ adX\circ adY)-Tr(adY\circ adH\circ adX)\}  \notag \\
&=&Tr\{adX\circ adH\circ adY)-(adX\circ AdY\circ adH)\}
\end{eqnarray}%
We have

\begin{equation}
\lbrack H,X]=2X,\text{ \ }[H,Y]=-2Y,\text{ \ }[X,Y]=H
\end{equation}%
Using (53) and (55), an easy computation gives $%
w_{3}^{G}(X,H,Y)=-8=-8v(X,H,Y)$ and therefore $w_{3}^{G}=-8\nu $. \ $\square 
$

Proposition 28 gives rise to a natural and, we believe, intriguing question.
Let $H\subset G$ be a subgroup of codimension $k.$ The cosets of $H$ foliate 
$G$ and this foliation descends to a foliation $\mathcal{F}$ on $G/\Delta $
with Godbillon-Vey class $GV(\mathcal{F})\in H^{2k+1}(G/\Delta ,\mathbb{R}).$
On the other hand we also have $[w_{2k+1}]_{dR}\in H^{2k+1}(G/\Delta ,%
\mathbb{R})$ which is defined even if $G$ does not have any subgroup of
codimension $k.$

$\mathbf{Q:}$ What is the relation between $GV(\mathcal{F})$ and $%
[w_{2k+1}]_{dR}$ ?

Since $[w_{1}]_{dR}$ is an obstruction to globalizability, it is natural to
expect that this is true for all $[w_{2k+1}]_{dR},$ that is, $%
[w_{2k+1}]_{dR}=0$ if $M$ is globalizable as we first conjectured. However,
the following example communicated to us by the referee shows that this
conjecture is false: let $G=SO(3,\mathbf{R})$ with the Lie algebra $%
\mathfrak{su}(3,\mathbb{R}\mathbf{)}$ generated by

\begin{equation}
A=\left[ 
\begin{array}{ccc}
0 & 1 & 0 \\ 
-1 & 0 & 0 \\ 
0 & 0 & 0%
\end{array}%
\right] ,\ B=\left[ 
\begin{array}{ccc}
0 & 0 & 1 \\ 
0 & 0 & 0 \\ 
-1 & 0 & 0%
\end{array}%
\right] ,\ C=\left[ 
\begin{array}{ccc}
0 & 0 & 0 \\ 
0 & 0 & 1 \\ 
0 & -1 & 0%
\end{array}%
\right]
\end{equation}

We have

\begin{equation}
\lbrack A,C]=B,\text{ }[B,A]=C,\ [C,B]=A
\end{equation}

As in the proof of Proposition 28, (57) together with (54 ) gives $%
w_{3}^{G}(A,B,C)=16\neq 0.$ Therefore $w_{3}^{G}$ is again a volume form and
is not exact since $G$ is compact.

The above examples $\mathfrak{sl}(2,\mathbb{R)}$, $\mathfrak{su}(3,\mathbb{R}%
\mathbf{)}$ are particular instances of a general phenomenon. To see this,
we continue to compute with the last equality in (62) by replacing $X,H,Y$
with arbitrary elements $x,y,z$ in some abstract Lie algebra $\mathfrak{g}$:

\begin{eqnarray}
w_{3}(x,y,z) &=&Tr\{adx\circ ady\circ adz)-(adx\circ adz\circ ady)\}  \notag
\\
&=&Tr\{adx\circ (ady\circ adz-adz\circ ady)  \notag \\
&=&Tr\{adx\circ ad[y,z]\}\text{ }  \notag \\
&=&\kappa (x,[y,z])
\end{eqnarray}

where $\kappa $ is Killing form.

\begin{proposition}
$i)$ $[w_{3}]_{\mathfrak{g}}\neq 0$ if $\mathfrak{g}$ is semisimple.

$ii)$ $w_{3}=0$ if and only if $\mathfrak{g}$ is solvable.
\end{proposition}

Proof: $i)$ This follows from (66) and the proof of Theorem 21.1 in [5].

$ii)$ This is again (66) together with the well known Cartan's criterion for
solvability. \ $\square $

The proof of Theorem 21.1 in [5] shows that $w_{3}$ is not only left
invariant but also right invariant, that is, invariant, and this is easy to
show for all $w_{2k+1}$ so that Proposition 27 is a consequence.

In view of Proposition 29, $[w_{3}]_{\mathfrak{g}}$ may be interpreted as an
obstruction to the solvability of the Lie algebra $\mathfrak{g}.$ Recalling
that the concept of solvability originated in Galois theory in an attempt to
solve polynomial equations by radicals, this interpretation is clearly far
from being satisfactory unless we clarify what is to be solved in the
present framework.

We would like to conclude with an amusing speculation about an intriguing
relation between local Lie groups and Poincare conjecture and a related
question. Let $M$ be a compact, simply connected 3-manifold together with
its unique differentiable structure. It is well known that $M$ is
parallelizable and therefore admits splittings by Lemma 3. If $M$ admits a
flat splitting, then it will be a local Lie group and therefore a Lie group
by Lemma 14 and Corollary 19, but $S^{3}$ is the only compact, simply
connected 3-dimensional Lie group. Therefore, Poincare Conjecture (now a
theorem thanks to the deep works of Hamilton and Perelman) is equivalent to
the assertion that $M$ admits a flat splitting. This fact suggests, we
believe, that Poincare Conjecture may have a short and metric-independent
proof.

The characteristic classes $[w_{2k+1}]_{dR}$ are secondary since they arise
when $M$ admits a splitting $\varepsilon $ with $\mathcal{R}(\varepsilon
)=0. $ The natural question arises whether there exist any primary classes
in the present framework.

$\mathbf{Q:}$ Are there any obstructions to the existence of a flat
splitting on a parallelizable manifold ?

\bigskip

\bigskip

\textbf{Acknowledgements: }We are indebted to the anonymous referee for
giving the counterexample $SO(3,\mathbf{R})$ to our conjecture that $%
[w_{2k+1}]_{dR}$ is an obstruction to globalizability also for $k\geq 1,$
and for pointing out that our "primary classes" vanish. We are also indebted
to Gregor Weingart for finding a serious mistake in the original version of
this paper and for his stimulating remarks. Last but not least, we express
our heartfelt gratitude to P. J. Olver for his kind support and
encouragement for this work.

\bigskip

\bigskip

\textbf{References}

\bigskip

[1] E. Abado\u{g}lu, E. Orta\c{c}gil, F. \"{O}zt\"{u}rk: Klein geometries,
parabolic geometries and differential equations of finite type, J. Lie
Theory 18 (2008) no.1, 67-82

[2] A.D. Blaom: Geometric structures as deformed infinitesimal symmetries,
Trans. Amer. Math. Soc. 358 (2006), no.8, 3651-3671(electronic)

[3] A.D. Blaom: Lie algebroids and Cartan's method of equivalence, arXiv:
math/0509071v2 [math.DG]

[4] R. Bott: Lectures on characteristic classes and foliations. Notes by
Lawrence Conlon, with two appendices by J.Stasheff. Lectures on algebraic
and differential topology (Second Latin American School in Math., Mexico
City, 1971, pp.1-94. Lecture Notes in Math., Vol. 279, Springer, Berlin, 1972

[5] C. Chevalley, S. Eilenberg: Cohomology theory of Lie groups and Lie
algebras, Trans. Amer. Math. Soc. 63 (1948), 85-124

[6] M. Crainic: Differentiable and algebroid cohomology, Van Est
isomorphisms, and characteristic classes, Comment. Math. Helv. 78 (2003),
no.4, 681-721

[7] M. Crainic, R.L. Fernandes: Secondary characteristic classes of Lie
algebroids, Quantum field theory and noncommutative geometry, 157-176,
Lecture Notes in Phys., 662, Springer, Berlin, 2005

[8] M. Crainic, I. Moerdijk: Deformations of Lie brackets: cohomological
aspects, arXiv:math / 0403434

[9] P.A. Damianou, R.L. Fernandes: Integrable hierarchies and the modular
class, Ann. Inst. Fourier (Grenoble) 58 (2008), no. 1, 107-137

[10] S. Evens, J.H. Lu, A. Weinstein: Transverse measures, the modular class
and a cohomology pairing for Lie algebroids, Quart. J. Math., Oxford Ser.
(2), 50 (1999), no. 200, 417-436

[11] R.L. Fernandes: Lie algebroids, holonomy and characteristic classes,
Adv. Math. 170 (2002), no. 1, 119-179

[12] R.B. Gardner: The Method of Equivalence and its Applications, SIAM,
Philedelphia, 1989

[13] P. Griffiths: On Cartan's method of Lie groups and moving frames as
applied to uniqueness and existence questions in differential geometry, Duke
Math. J.. 41 (1974), 775-814

[14] S. Helgason: Differential geometry, Lie groups, and symmetric spaces.
Corrected reprint of the 1978 original, Graduate Studies in Mathematics, 34.
AMS, Providence, RI, 2001

[15] F.W. Kamber, P. Tondeur: Foliated bundles and characteristic classes.
Lecture Notes in Mathematics, Vol. 493. Springer-Verlag, Berlin-New York,
1975

[16] S. Kobayashi, K. Nomizu: Foundations of Differential Geometry. Vol. I,
Interscience Publishers, a division of John Wiley \& Sons, New York-London,
1963

[17] A. Kumpera, D.C. Spencer: Lie Equations. Vol. I: General Theory, Annals
of Mathematics Studies, No. 73, Princeton University Press, Princeton, N.J,
Univertsity of Tokyo Press, Tokyo, 1972

[18] K. MacKenzie: Lie Groupoids and Lie Algebroids in Differential
Geometry, London Mathematical Society Lecture Note Series, 124, Cambridge
University Press, Cambridge, 1987

[19] K. MacKenzie: General Theory of Lie Groupoids and Lie Algebroids,
London Mathematical Society Lecture Note Series, 213, Cambridge University
Press, Cambridge, 2005

[20] P.J. Olver: Non-associative local Lie groups, J. Lie \ 6 (1996), no.1,
23-51

[21] E. Ortacgil: The heritage of S.Lie and F.Klein: Geometry via
transformation groups, arXiv:math / 0604223, April 2006

[22] E. Orta\c{c}gil: Pre-homogeneous geometric structures and their
curvatures, in progress

[23] J.F. Pommaret: Systems of Partial Differential Equations and Lie Pseudo

groups. With a preface by Andre Lichnerowicz, Mathemetics and its
Applications, 14. Gordon\&Breach Science Publishers, New York, 1978

[24] J.F. Pommaret: Partial Differential Equations and Group Theory. New
Perspectives for Applications, Mathematics and its Applications, 293. Kluwer
Academic Publishers Group, Dordrecht, 1994

[25] L.S. Pontryagin: Selected works. Vol.2, Topological groups, Edited and
with a preface by R.V. Gamkrelidze, translated by Arlen Brown007 (90a:01106)
Pontryagin, L. S. Selected works. Vol. 2. Topological groups, with a preface
by R. V. Gamkrelidze, Classics of Soviet Mathematics, Gordon \& Breach
Science Publishers, New York, 1986

[26] R.W. Sharpe: Differential Geometry, Cartan's Generalization of Klein's
Erlangen Program, with a foreword by S.S.Chern, Graduate Texts in
Mathematics, 166. Springer-Verlag, New York, 1997

[27] O. Veblen: Invariants of quadratic differential forms, Cambridge
University Press, 1962

[28] G. Weingart: Holonomic and semi-holonomic geometries, Global Analysis
and Harmonic Analysis (Marseille-Luminy, 1999), 307-328, Semin. Congr., 4,
Soc. Math. France, Paris, 2000

[29] A. Weinstein: The modular automorphism group of a Poisson manifold, J.
Geom. Phys., 23 (1997), no. 3-4, 379-394

\bigskip

Ender Abado\u{g}lu, Yeditepe University, Mathematics Department, 81120,
Kay\i \c{s}da\u{g}\i , \.{I}stanbul, T\"{u}rkiye

e-mail: eabadoglu@yeditepe.edu.tr

\bigskip

Erc\"{u}ment Orta\c{c}gil, Bo\u{g}azi\c{c}i University, Mathematics
Department, Bebek, 34342, \.{I}stanbul, T\"{u}rkiye

e-mail: ortacgil@boun.edu.tr

\end{document}